\documentclass[12pt]{article}
\usepackage{amscd,amsfonts,amssymb,amscd,amsmath,latexsym}
\usepackage[hypertex,backref]{hyperref}
\makeatletter
\renewcommand{\subsubsection}{\@startsection
{subsubsection}
{1}
{0mm}
{0mm}
{0mm}
{\normalfont\normalsize\itshape}}
\makeatother
  
\usepackage{mathrsfs}
\newtheorem{theorem}{Theorem}[section] 
\newtheorem{prop}[theorem]{Proposition}
\newtheorem{lem}[theorem]{Lemma}
\newtheorem{ddd}[theorem]{Definition}
\newtheorem{kor}[theorem]{Corollary}

\newcommand{\Lie}{{\tt Lie}}

\newcommand{\forget}[1]{}

{\catcode`@=11\global\let\c@equation=\c@theorem}

\newcommand{\Z}{\mathbb{Z}}

\newcommand{\proof}{{\it Proof.$\:\:\:\:$}}

\newcommand{\R}{\mathbb{R}}

\newcommand{\Cliff}{{\tt Cliff}}

\newcommand{\C}{\mathbb{C}}

\newcommand{\Aut}{{\tt Aut}}

\newcommand{\Tr}{{\tt Tr}}
\newcommand{\cT}{\mathcal{T}}

\newcommand{\cH}{\mathcal{H}}
\newcommand{\cM}{\mathcal{M}}

\newcommand{\hol}{{\tt hol}}

\newcommand{\cI}{\mathcal{I}}

\newcommand{\Hom}{{\tt Hom}}

\newcommand{\Ext}{{\tt Ext}}

\newcommand{\sign}{{\tt sign}}

\newcommand{\cN}{\mathcal{N}}

\newcommand{\Ad}{{\tt  Ad}}

\newcommand{\codim}{{\tt codim}}

\newcommand{\id}{{\tt id}}

\def\imath{{i}}
\newcommand{\cR}{\mathcal{R}}
\def\hB{\hspace*{\fill}$\Box$ \newline\noindent}

\newcommand{\cS}{\mathcal{S}}

\def\hB{\hspace*{\fill}$\Box$ \\[0cm]\noindent}
\newcommand{\cL}{\mathcal{L}}
 
\newcommand{\cP}{\mathcal{P}}
\newcommand{\cQ}{\mathcal{Q}}

\newcommand{\pr}{{\tt pr}}

\title{Twisted $K$-theory and TQFT}
\author{Ulrich Bunke and Ingo Schr{\"o}der
\thanks{Mathematisches Institut, Universit{\"a}t G{\"o}ttingen,
Bunsenstr. 3-5, 37073 G{\"o}ttingen, GERMANY,  bunke@uni-math.gwdg.de,
ischroed@uni-math.gwdg.de}
}

\begin{document}
 

\maketitle 
\tableofcontents
\parskip3ex

\section{Introduction}

The present paper grew out of a seminar held in spring, 2004. The goal
of the seminar was to understand the recent paper by Freed, Hopkins,
and Teleman \cite{freedhopkinsteleman03}. The first main result  of
\cite{freedhopkinsteleman03} that we
discussed was the calculation of the twisted $G$-equivariant $K$-theory
of $G$, where $G$ is a compact Lie group which  acts on itself by conjugation. 

While working on details we
came to the conclusion that it is  worth to develop proofs in a more
restricted formalism. The original paper \cite{freedhopkinsteleman03}
mixes analytic with geometric and topological arguments. We had some
difficulties to see that all constructions match in a nice manner.

In the present paper we try to give a proof of this result (which we
formulate here as Theorem
\ref{theo:calulation}) by  arguments which are completely embedded
in the calculus of smooth stacks. We do not touch the question of the
construction of a $K$-theory functor in this framework. Rather we
assume that such a functor exists and has  all necessary functorial
properties. Actually we only need local quotient stacks, and the
construction of the $K$-theory in this case was sketched in
\cite{freedhopkinsteleman03} (see also \cite{laurenttuxu04} and \cite{atiyahsegal04}). 
A verification of all functorial properties, in
particular the construction of push-forward maps, is still a gap in
the literature. 

The way the calculation is set up in the present paper opens the path to
generalizations. Since we only do geometric calculations with stacks 
and use the formal properties of $K$-theory, the method could also be
applied to other twisted cohomology theories.

The calculation of the twisted $K$-theory has two basis steps. The 
first step is the   construction of  elements of the twisted
$K$-theory. In symbols the construction is realized as the map
$R_!\circ \Phi$, where $\Phi$ is defined in \ref{phidef41}, and the
map $R$ is introduced in \ref{restor}. Note that this construction is
purely geometric in terms of the calculus of smooth stacks and formal
properties of twisted $K$-theory. 

The second step is a method to detect elements of the twisted
$K$-theory. We will construct an embedding (the map $\Theta$
introduced in \ref{thetadef}, see Theorem \ref{detec})  of the twisted $K$-theory
into the representation ring of a suitable finite group \ref{ctd1}.
This very effective tool was explained to us by C. Teleman.

The second topic of the seminar were elements of a $1+1$-dimensional TQFT-structure on the twisted
$K$-theory. The identity, the product, and the co-form are induced by
natural geometric constructions with stacks associated to the group.
Having constructed a basis of the twisted $K$-theory in Theorem
\ref{theo:calulation} it is then a natural question to express these
TQFT-operations in terms of this basis. The results and sketches of
proofs were announced again in \cite{freedhopkinsteleman03} and
\cite{fht04}. In the present paper we reproduce the formulas working
again completely inside the stack calculus and using only formal
properties of $K$-theory.
The identity is calculated in Theorem \ref{theunit}. The co-form is
obtained in Theorem \ref{thedi2}. Finally, the calculation of the product is stated as
Theorem \ref{pru78}.

The twisted $K$-theory of a Lie group which acts on itself by
conjugation is a module over the representation ring of the Lie group
in a natural way. It follows from the calculation that the twisted
$K$-theory is a free  $\Z$-module and therefore
embeds in its complexification. The latter is a module over the
complexified group ring. In Theorem \ref{rrg1} we show that it is a
quotient of the complexified group ring. Actually, equipped with the
product and the identity coming from the TQFT-structure, the
complexified twisted $K$-theory is a quotient of the complexified
representation ring as a ring.

The natural source of the TQFT-structure  are
correspondences given by moduli spaces
of flat connections on surfaces and their boundaries (see
(\ref{cne3})). In the present paper these correspondences were only
employed to construct isomorphisms of twists needed to define the product.

Indeed, we had difficulties to define the $K$-orientations of the
outgoing boundary evaluation maps in a natural way such that they are
compatible with glueing.

In Subsection \ref{pol91} we recall the construction of the central extension of the
restricted unitary group associated to a polarized Hilbert spaces.
This central extension can be viewed as a source of twists. In
\ref{uf49} we use this central extension in order define in a natural way  twists of the moduli stack of (flat)
$G$-connections on any one-dimensional compact closed oriented manifold.

The determinant line bundle over the restricted Grassmannian of the
polarized Hilbert space can be used as a source of trivializations of
twists. For a compact oriented surface with boundary we have an
evaluation map from the moduli stack of flat $G$-connections on the
surface to the stack of $G$-connections on the boundary.
In  \ref{mmc1} we construct a natural
trivialization of the pull-back of the twist via the boundary
evaluation. In Propositions \ref{suin} and  \ref{rrc17}
we verify that these trivializations behave functorially with respect
to the glueing of surfaces. This approach to twists was partly inspired by the thesis of Posthuma \cite{posthuma}. 

The missing piece for a completely natural
construction of the TQFT using moduli spaces is the compatible
orientation of the outgoing boundary evaluation maps.
Such a construction is desirable in particular, because it would give a natural explanation for the
associativity of the product. 

A major topic of \cite{freedhopkinsteleman03} is the relation between the
equivariant twisted $K$-theory of the Lie group acting on itself by
conjugation and the theory of positive energy representations of the associated
loop group. Because of lack of time this was not discussed in the
seminar and will therefore not be touched in the present
paper. Another more philosophical reason for this omission is that according to our
present knowledge this relation
can not be seen purely inside the calculus of stacks. Rather it is based
on explicite cycles in order to represent twisted $K$-theory classes
in an appropriate model.

While working on this paper we had a very fruitful exchange with
C. Teleman. He told us the idea how to detect elements of the twisted
$K$-theory groups using the restriction to finite groups. Furthermore 
this discussion led to the elimination of many stupid mistakes in
previous versions of these notes.

In the same seminar J. Heinloth gave an introduction to smooth stacks
and gerbes. In the present paper we freely use the language and the
notation which was set up in his talks and the review \cite{heinloth}. Further
discussions with J. Heinloth during the preparation of the present
paper were of great help.

Finally, since this will not be noted again in the text below, let us emphasize that the main theorems discussed in
the present paper and the key ideas leading to their verifications are due to
Freed, Hopkins, and Teleman.

\section{Calculation of twisted $K$-theory of Lie groups}

\subsection{Connections, gauge groups, and twists}

\subsubsection{}\label{str43}

Let $G$ be a Lie group. We consider the trivial $G$-principal bundle
$$P(S^1):=G\times S^1\rightarrow S^1\ .$$ Let $F(S^1)$ denote the space of
(flat) connections on $P(S^1)$. The gauge group $G(S^1)$ acts on $F(S^1)$.
We consider the topological stack (see \cite{heinloth}, Ex. 1.5 and 2.5) $$\cM:=[F(S^1)/G(S^1)]\ .$$

\subsubsection{}

Let $$0\rightarrow U(1)\rightarrow \hat G(S^1)\rightarrow
G(S^1)\rightarrow 0$$ be a central extension. It gives rise to a twist
(see \cite{heinloth}, Ex. 5.4.1)
$$\tau:\hat \cM\rightarrow \cM\ ,$$
where $\hat \cM:=[F(S^1)/\hat G(S^1)]$.
The goal of the present section is a to formulate the main result
about the calculation of the twisted
$K$-theory ${}^\tau K(\cM)$.

\subsubsection{}

Twisted $K$-theory associates to a topological stack $\cM$ equipped
with a twist $\tau:\hat\cM\rightarrow \cM$ a $\Z$-graded group
${}^\tau K(\cM)$ in a functorial way. More precisely, if
$\tau^\prime:\hat \cM^\prime\rightarrow \cM^\prime$ is another twisted
topological stack, $f:\cM^\prime\rightarrow \cM$ is a morphism, and
$u:\tau^\prime\rightarrow f^*\tau$ is an isomorphism of twists, then we
have a functorial map $u^*f^*:{}^\tau K(\cM)\rightarrow
{}^{\tau^\prime}K(\cM^\prime)$. See
\cite{freedhopkinsteleman03}, \cite{laurenttuxu04}, and
\cite{atiyahsegal04} for a construction a twisted $K$-theory functor.
We further assume that twisted $K$-theory admits a Mayer-Vietoris
sequence and is a module over the untwisted $K$-theory.
Our assumptions about wrong-way maps will be  explained in \ref{asor5}.

\subsubsection{}\label{somp6}

We consider $G\subset G(S^1)$ as the subgroup of constant gauge transformations.
We assume that $G$ is connected and
choose a maximal torus $T\subset G$. Let $\check{T}$ denote the group
of homomorphisms $S^1\rightarrow T$. We can consider $\check{T}\subset
G(S^1)$ naturally. Furthermore let $N_G(T)$ be the normalizer of $T$
in $G$ which we also consider as a subgroup of $G(S^1)$. 
Inside $G(S^1)$ the groups $\check{T}$ and $N_G(T)$ generate a
semi-direct product 
$$0\rightarrow \check{T}\rightarrow \check{T}N_G(T)\rightarrow
N_G(T)\rightarrow 0\ .$$
The group of connected components of $
\check{T}N_G(T)$ is the affine Weyl group $\hat W$.
It fits into a semi-direct product
$$0\rightarrow \check{T}\rightarrow \hat W\rightarrow
W\rightarrow 0\ ,$$
where $W:=N_G(T)/T$ is the ordinary Weyl group of the pair $(G,T)$.

\subsubsection{}

Let $\hat T\rightarrow T$ be the restriction of the central extension
of $G(S^1)$ via the embedding $T\subset G(S^1)$. Let $X(\hat T)$
denote the group of characters, and let $X_1(\hat T)\subset X(\hat T)$ be
the subset of those characters which become the identity after
restriction to the central $U(1)$.

The torus $T\subset G(S^1)$ is preserved under conjugation by elements
of $\check{T}N_G(T)$. Therefore $\check{T}N_G(T)$ acts on $X(\hat T)$.
In fact, this action preserves $X_1(\hat T)$ and factors over the
affine Weyl group.

\subsubsection{}\label{reg5412}

We call an element $\chi\in X_1(\hat T)$ regular, if its
stabilizer in $\hat W$ is trivial. Otherwise we call $\chi$ singular.
Let $X^{reg}_1(\hat T)$ denote the set of regular elements.

We call the twist $\tau$ regular, if $\hat W$ acts properly on
$X_1(\hat T)$ with finitely many orbits. 

Let $\widehat{\check{T}N_G(T)}\rightarrow 
\check{T}N_G(T)$ be the central extension induced by the restriction of
$\hat G(S^1)\rightarrow G(S^1)$ to $\check{T}N_G(T)$. By further restrictions we obtain
central extensions $\widehat{N_G(T)}$ and $\widehat{\check{T}}$ of
$N_G(T)$ and $\check{T}$.
We call the twist $\tau$ admissible if $\widehat{N_G(T)}$ and $\widehat{\check{T}}$ are trivial.

\subsubsection{}

The main result of the
present section is the formulation of the following theorem:
\begin{theorem}\label{theo:calulation} Assume that $G$ is connected, and that $\tau$ is
  regular and admissible. Then
the orbit  set $X_1^{reg}(\hat T)/\hat W$ is the index set of a $\Z$-basis
of the free $\Z$-module ${}^\tau K(\cM)$ in a natural way.
\end{theorem}
We will finish the proof of this theorem in  \ref{thebases}.
After a choice of representatives of the equivalence classes  
$X_1^{reg}(\hat T)/\hat W$ the basis elements will be determined uniquely up to a global
sign which can be fixed by choosing an orientation of $\Lie(T)$.

\subsubsection{}

Let $[G/G]$ be the quotient stack, where $G$ acts in itself by conjugation.
We define a map $\hol:\cM\rightarrow [G/G]$ which on the level of
spaces associates to each connection in $F(S^1)$ its holonomy at $1\in
S^1$ measured in the positive direction. On the level of groups it is
given by the evaluation $G(S^1)\rightarrow G$ at $1$.

\begin{lem}
The map $\hol: \cM\rightarrow [G/G]$ is an equivalence of topological stacks.
\end{lem}
\proof
This follows from \cite{heinloth}, Ex. 3.3 and the fact that the group
$G(S^1)_0\subset G(S^1)$ of based gauge transformations (those which
evaluate trivially at $1\in S^1$) acts freely and properly on $F(S^1)$
with quotient isomorphic to $G$ via the holonomy map, and
$G(S^1)/G_0(S^1)\cong G$ via the evaluation.
\hB

The composition $\hol_*\tau:=\hol\circ \tau:\hat \cM\rightarrow [G/G]$
is a twist of $[G/G]$.
In this sense Theorem \ref{theo:calulation} provides a calculation of
${}^{\hol_*\tau} K( [G/G])$.

\subsubsection{}\label{oddef1}

We consider a regular and admissible twist $\tau$ of $\cM$. Let $I:[G/G]\rightarrow
[G/G]$ be the map which is given by $g\mapsto g^{-1}$ on the level of
spaces, and by the identity on the level of groups.
We call the twist $\tau$  odd, if $I^*\hol_*\tau\cong -\hol_*\tau$.


\subsection{Orientations}\label{pokl}

\subsubsection{}

Let $E$  be a real euclidian vector space. By $\Cliff(E)$ we denote
the  associated complex Clifford algebra. It comes with an embedding
of $E\rightarrow \Cliff(E)$ and a $*$-operation.
Let $\Cliff(E)^*$ denote the group of unitary elements.
We define
$$Pin^c(E):=\{x\in\Cliff(E)^*\:|\: xEx^{*}=E\}\ .$$
This group comes as a central extension
\begin{equation}\label{seqgr}0\rightarrow U(1)\rightarrow Pin^c(E)\rightarrow O(E)\rightarrow 0\
.\end{equation}
If $E=\R^n$, then we write $Pin^c(n):=Pin^c(\R^n)$.

We let $Spin^c(E)\subset Pin^c(E)$ be the preimage of $SO(E)\subset O(E)$, and set
$Spin^c(n):=Spin^c(\R^n)$.

\subsubsection{}

The sequence of groups (\ref{seqgr}) induces a sequence of maps of stacks
$$ 
[*/U(1)]\rightarrow [*/Pin^c(n)]\rightarrow [*/O(n)]\\
 \ .$$
We furthermore have the following pull-backs (see \cite{heinloth},
Def. 2.1)
$$\begin{array}{ccccc}
{}[*/Spin^c(n)]&\rightarrow&[*/SO(n)]&\rightarrow&*\\
\downarrow&&\downarrow&&\downarrow\\
{}[*/Pin^c(n)]&\rightarrow&[*/O(n)]&\rightarrow&[*/(\Z/2\Z)]\end{array}\
.$$

\subsubsection{}\label{twex}

Let $E\rightarrow \cM$ be a real euclidian vector bundle over a stack
(see \cite{heinloth}, 2.10).
Its frame bundle gives rise to a classifying map
$\cM\rightarrow[*/O(n)]$. We form the pull-back
$$\begin{array}{ccc}
Spin^c(E)&\rightarrow&[*/Spin^c(n)]\\
\downarrow&&\downarrow\\
\cM&\rightarrow&[*/O(n)]\end{array}\ .$$
The stack $Spin^c(E)$ classifies $Spin^c$-structures on $E$.
We compose the map $\cM\rightarrow [*/O(n)]$ with
$[*/O(n)]\rightarrow [*/(\Z/2\Z)]$ and obtain
the pull-back 
$$\begin{array}{ccc}
or(E)&\rightarrow&*\\
\downarrow&&\downarrow\\
\cM&\rightarrow &[*/(\Z/2\Z)]
\end{array}$$
defining the stack $or(E)$ which classifies the
orientations of $E$. Then we have a natural factorization
$Spin^c(E)\rightarrow or(E)$ which is a  twist. Below will call the
pair $(or(E)\rightarrow \cM ,Spin^c(E)\rightarrow or(E))$ a graded
twist and denote it by $\tau(E)$.

\subsubsection{}

A graded twist of a stack $\cM$ will be a pair
$(\sigma,\rho)$, where $\sigma:\cM\rightarrow [*/(\Z/2\Z)]$ (is called
the grading),
and $\rho$ is a twist
of the $\Z/2\Z$-covering $\sigma(\cM)$ defined by
$$\begin{array}{ccc}
\sigma(\cM)&\rightarrow&*\\
\downarrow&&\downarrow\\
\cM&\stackrel{\sigma}{\rightarrow} &[*/(\Z/2\Z)]
\end{array}\ .$$

\subsubsection{}\label{taud}

A real Euclidean vector bundle $E\rightarrow \cM$ gives rise to a
graded twist
$\tau(E)$ as explained in \ref{twex}.
The twists of a stack form
a two-category. 

There is a natural notion of a sum of  graded
twists such that there is a natural isomorphism
$\tau(E_0\oplus E_1)\cong \tau(E_0)+\tau(E_1)$.
This sum differs from the componentwise sum. In fact, the isomorphism
classes of twists are classified by a group $GTW(\cM)$ which sits in an non-trivial
extension
$$0\rightarrow H^3(\cM,\Z)\rightarrow GTW(\cM)\rightarrow
H^1(\cM,\Z/2\Z)\rightarrow 0$$
(see the discussion in \cite{atiyahsegal04} and \cite{heinloth},
Remark 5.9).

\subsubsection{}\label{cohe}

 Let $f:\cM\rightarrow \cN$ be a representable smooth map of smooth
stacks. A factorization of $f$ into 
a smooth embedding and a smooth
submersions gives rise to a normal bundle (see \cite{heinloth}, Def. 4.10)
which is a $\Z/2\Z$-graded vector bundle $E=E^+\oplus E^-$ over $\cM$.
We define $\tau(E):= \tau(E^+)-\tau(E^-)$.
Let $E^\prime$ be a normal bundle obtained by a different
faktorization. Then using the diagonal embedding we get  bundles $A,B$ over $\cM$ such that
$E\oplus A \oplus (-A)$ and $E^\prime\oplus B\oplus (-B)$ are
canonically isomorphic. In particular, we get a
natural isomorphism $\tau(E)\cong \tau(E^\prime)$.




\subsubsection{}

Let $f:\cM\rightarrow \cN$ be a representable smooth map which admits
factorizations in smooth embeddings and submersions. Let $\tau$ be a
graded twist of $\cN$.

\begin{ddd}
A $\tau$-$K$-orientation of $f$ is a coherent choice of
isomorphisms $$f^*\tau\stackrel{\sim}{\rightarrow} \tau(E)$$
of twists for all normal bundles $E$ of $f$ given by some
factorization.
\end{ddd}

Coherence is understood here with respect to the natural isomorphisms
(see \ref{cohe}) 
of twists associated to the normal bundles given by different
factorzations of $f$.
Note that a $\tau$-$K$-orientation is determined by the isomorphism
$f^*\tau\stackrel{\sim}{\rightarrow} \tau(E)$ for one choice of $E$.

\subsubsection{}

Let $$\begin{array}{ccc}
\cQ&\stackrel{q}{\rightarrow}&\cN\\
g\downarrow&&f\downarrow\\
\cR&\stackrel{p}{\rightarrow}&\cM\end{array}$$
be a cartesian diagram of smooth stacks, where $f$ is representable and has a
factorization with normal bundle $E$. Then $g$ is representable and
also has a factorization with normal
bundle $q^*E$. If $f$ is $\tau$-$K$-oriented, then
$g$ has an induced $p^*\tau$-$K$-orientation.
 
\subsubsection{}\label{twoofthree}

Let us consider a composition
$$\cR\stackrel{g}{\rightarrow}\cN\stackrel{f}\rightarrow \cM\ .$$
If $E$ and $F$ are normal bundles (for some factorizations) of $f$ and
$g$, then there is a factorization  of $f\circ g$  with a normal
bundle  $g^*E\oplus F$. 
If $f$ is $\tau_0$-$K$-oriented, and $g$ is
$f^*\tau_1$-oriented, then
$f\circ g$ is naturally $ \tau_0+\tau_1$
oriented.  Vice versa, if $f$ and $f\circ g$ are $K$-oriented, then so is
canonically $g$.

\subsubsection{}\label{asor5}

We assume that the twisted $K$-theory functor admits functorial wrong-way maps for
twisted-$K$-oriented proper maps. More precisely, if $f:\cM\rightarrow \cN$
is a smooth proper (see \cite{heinloth}, Def. 2.8) map between stacks
which is $\tau$-$K$-oriented for some twist $\tau$ of $\cN$,  and if
$\sigma$ is a further twist of $\cN$, then we have a wrong-way map
$$f_!:{}^{f^*\sigma+f^*\tau}K(\cM)\rightarrow {}^{\sigma} K(\cN)\ .$$
We assume functoriality with respect to compositions and compatibility
with cartesian diagrams. Furthermore, we require a projection formula 
for the module structure of the twisted $K$-theory over the untwisted $K$-theory.

\subsubsection{}\label{orfix}

Let $G$ be a connected compact Lie group. We consider $G$ as a
$G\times G$-space with with the action
$(a,b)g=agb^{-1}$.

The group $G$ acts on its Lie algebra $\Lie(G)$ by the adjoint
representation. We let $G\times G$ act on $\Lie(G)$ via the projection
onto the second factor.

We obtain a vector bundle $[\Lie(G)/G\times G]\rightarrow [*/G\times
G]$. Let $\tau(G)$ be the corresponding obstruction twist as in
\ref{taud}. In fact, since $G$ is connected we do not need the grading.

We will fix once and for all an orientation of $\Lie(G)$. It induces
an orientation of $G$. Since $G$ is connected the orientation covering
$or(\Lie(G))\rightarrow [*/G\times G]$ is trivialized.


\subsubsection{}
\begin{lem}\label{kors}
The map $q:[G/G\times G ]\rightarrow [*/G\times G]$ has a natural  $-\tau(G)$-$K$-orientation.
\end{lem}
\proof
The normal bundle of $q$ is the tangent bundle
$[TG/G]\rightarrow[G/G]$ placed in degree one.
We trivialize the
tangent bundle $TG\cong G\times \Lie(G)$ using the left action. The action of $(g_1,g_2)\in
G\times G$ on $(h,X)\in  G\times \Lie(G)$ is then given by
$(g_1hg_2^{-1},\Ad(g_2)(X))$. Therefore we obtain an isomorphism
$[TG/G\times G]\cong q^*[\Lie(G)/G\times G]$. This gives the natural
$-\tau(G)$-$K$-orientation of $p$.   \hB



\subsubsection{}\label{my1}

Let $d:[*/G]\rightarrow [*/G\times G]$
be given  by the diagonal embedding $G\rightarrow G\times G$.
Then we have a cartesian diagram
$$\begin{array}{ccc}
[G/G]&\rightarrow&[G/G\times G]\\
p\downarrow&&q\downarrow\\{}[*/G]&\stackrel{d}{\rightarrow}&[*/G\times
G]\end{array}\ .
$$
The $-\tau(G)$-$K$-orientation of $q$ induces a
  $-\sigma(G)$-$K$-orientation of $p:[G/G]\rightarrow [*/G]$,
where $\sigma(G):=d^*\tau(G)$.
 We consider the sequence
$$[*/G]\stackrel{i}{\rightarrow} [G/G]\stackrel{p}{\rightarrow}[*/G]\ .$$
The canonical $K$-orientation of the composition $p\circ i=\id$ and the
$-\sigma(G)$-$K$-orientation of $p$ induce a $p^*\sigma(G)$-$K$-orientation of $i$.

\subsubsection{}\label{my2}

The normal bundle of the map
$*\rightarrow [*/G]$ is 
$\Lie(G)\rightarrow *$ placed in degree one. We have already fixed an orientation in \ref{orfix}.
The unique $Spin^c$-structure on the vector bundle $\Lie(G)\rightarrow *$
induces the $K$-orientation of $$*\rightarrow [*/G]\ .$$




\subsubsection{}\label{sigmatr12}

Twists of $[*/G]$ are classified by $H^3([*/G],\Z)\cong H^3(BG,\Z)$
(see \cite{heinloth}, Prop. 5.8).
In fact, the class of $\sigma(G)$ is two-torsion since it comes from a
finite-dimensional vector bundle. Note that
$H^3(BG,\Z)_{tors}\cong \Ext(H_2(BG,\Z),\Z)$, and that $H_2(BG,\Z)\cong
\pi_1(G)$. Therefore, if we assume that $2$ does not divide the order
of $\pi_1(G)_{tors}$, then
$\sigma(G)$ is trivial.

Note further, that the isomorphism classes of trivializations of
$\sigma(G)$ form a $H^2(BG,\Z)$-torsor (see \cite{heinloth}, Remark 5.7), and that $H^2(BG,\Z)\cong
\Hom(\pi_1(\Z),\Z)$ (since $G$ is connected). Thus, if we assume that
$\pi_1(G)$ is finite, then $\sigma(G)$ is trivial in a unique way.

\subsubsection{}

Let $T\subset G$ be a maximal torus and $N_G(T)$ be its normalizer.
\begin{lem}\label{tori}
The map $q:[T/N_G(T)]\rightarrow [*/G]$ has a natural $-\sigma(G)$-$K$-orientation.
\end{lem}
\proof
We represent (see \cite{heinloth}, Ex. 3.3) this map as
$$[(G\times_{N_G(T)}T)/G]\rightarrow [*/G]\ .$$
Then we can write the tangent bundle of  $G\times_{N_G(T)}T$ as associated
vector bundle
$$(G\times T)\times_{N_G(T)} \{(\Lie(G)/\Lie(T))\oplus \Lie(T)\}\ .$$
Since the representation of $N_G(T)$ on $\Lie(G)/\Lie(T)\oplus \Lie(T)\cong
\Lie(G)$ extends to the adjoint representation of $G$ we have an
isomorphism of $G$-equivariant bundles
$$(G\times T)\times_{N_G(T)} \{(\Lie(G)/\Lie(T))\oplus \Lie(T)\}\cong
(G\times_{N_G(T)} T)\times \Lie(G)\ .$$
Therefore we can identify the vertical bundle of $q$ with the
pull-back by $q$
of $\Lie(G)\rightarrow *$ as $G$-equivariant bundles. This provides
the natural $-\sigma(G)$-$K$-orientation of $q$.
\hB

\subsubsection{}\label{orir}

We consider the composition
$$[T/N_G(T)]\stackrel{R}{\rightarrow} [G/G]\stackrel{\pi}{\rightarrow}
[*/G]\ ,$$
where $R$ is induced by the obvious embeddings on the level of spaces
and groups.
Now $\pi$  has a natural $-\sigma(G)$-$K$-orientation by Lemma \ref{kors}, and the
composition $\pi\circ R$ has a natural $-\sigma(G)$-$K$-orientation by
Lemma \ref{tori}. It follows that the map $R$ has a natural $K$-orientation.

\subsection{Construction of twisted $K$-theory}\label{restor}

\subsubsection{}

We have an embedding $\Lie(T)\rightarrow F(S^1)$ as constant
connections. The group $\check{T}N_G(T)\subset G(S^1)$ preserves the
image and thus acts on $\Lie(T)$. Let
$\cT:=[\Lie(T)/\check{T}N_G(T)]$.
Then we have a map of stacks $$R:\cT\rightarrow \cM\ .$$

Let $\widehat{\check{T}N_G(T)}\rightarrow \check{T}N_G(T)$ denote the
restriction of the central extension $\hat G(S^1)\rightarrow G(S^1)$ and set $$\hat
\cT:=[\Lie(T)/\widehat{\check{T}N_G(T)}]\ .$$ Then we have the twist
(see \cite{heinloth}, Remark 5.5.3 for the pull-back of a twist)
$$R^*\tau:\hat \cT\rightarrow \cT\ .$$



\subsubsection{}

We consider the $\hat T$-principal bundle
$h:\widehat{\check{T}N_G(T)}\rightarrow \hat W$, where $h$ is the
projection to the group of connected components. We let $L^2(h)\rightarrow \hat W$ be the bundle of Hilbert spaces
such that its fibre $L^2(h)_{\hat w}$ over $\hat w\in \hat W$ is 
$L^2(h^{-1}(\hat w))$.
Note that $\hat T$ acts on $L^2(h)_{\hat w}$ via the right action on the fibre.
  We further define the line bundle
$L\rightarrow \hat W\times X_1(\hat T)$ such that
the fibre $L_{(\hat w,\chi)}$ over $(\hat w,\chi)$ is the $\chi$-isotypic component 
$L^2(h)_{\hat w}(\chi)\subset L^2(h)_{\hat w}$ as a $\hat T$-representation. 



We define an action of $$\widehat{\check{T}N_G(T)}\times \widehat{\check{T}N_G(T)}$$ on $$L\rightarrow \hat W\times X_1(\hat T)$$ as follows.
Let 
$f\in L^2(h)_{\hat w}(\chi)=L_{(\hat w,\chi)}$. Then we define
$(\hat n_1,\hat n_2)(f)(\hat g):= f(\hat n_1^{-1} \hat g \hat n_2)$ for all $\hat g\in h^{-1}(\hat n_1\hat w\hat n_2^{-1})$.
We calculate that
$(\hat n_1,\hat n_2)(f)(\hat g \hat t)= f(\hat n_1^{-1} \hat g \hat t\hat n_2)=
f(\hat n_1^{-1} \hat g \hat n_2 \hat n_2^{-1}\hat t\hat n_2)=f(\hat n_1^{-1} \hat g \hat n_2)\chi(\hat n_2^{-1}\hat t\hat n_2)=
(\hat n_2\chi)(\hat t) (\hat n_1,\hat n_2)(f) (\hat g) 
$. It follows that
$(\hat n_1,\hat n_2)(f)\in L_{(\hat n_1\hat w\hat n_2^{-1},\hat n_2\chi)}$.
Therefore the projection of $L$ becomes equivariant if we let  
$\widehat{\check{T}N_G(T)}\times \widehat{\check{T}N_G(T)}$ act on the
base $\hat W\times X_1(\hat T)$
by $(\hat n_1,\hat n_2)(\hat w,\chi):=(\hat n_1\hat w\hat n_2^{-1},\hat n_2\chi)$. Since we assume that $\tau$ is admissible we can choose a split
$\hat W\rightarrow \widehat{\check{T}N_G(T)}$ of $h$ which is a homomorphism.
We let $\hat W$ act on $L\rightarrow \hat W\times X_1(\hat T)$ via its
embedding into the left factor of $\widehat{\check{T}N_G(T)}\times
\widehat{\check{T}N_G(T)}$ given by the split. It acts freely, and the quotient is a certain
$\widehat{\check{T}N_G(T)}$-equivariant line bundle 
$$\bar L\rightarrow X_1(\hat T)\ .$$ Note that the central $U(1)\subset
\widehat{\check{T}N_G(T)}$ acts on the fibres of $\bar L$ by the
identity character ($\bar L$ is of weight one in the language of
\cite{heinloth}, Lemma 5.6)
Therefore we can consider the unit sphere bundle
$$[U(\bar L)/\widehat{\check{T}N_G(T)}]\rightarrow [X_1(\hat T)/\widehat{\check{T}N_G(T)}]$$ as a trivialzation of the twist
$$[X_1(\hat T)/\widehat{\check{T}N_G(T)}]\rightarrow [X_1(\hat
T)/\check{T}N_G(T)]\ .$$


\subsubsection{}\label{gg23}

We consider the space
$\Lie(T)\times X_1(\hat T)$ with the diagonal action of
$\check{T}N_G(T)$.
It gives rise to the stack $$\cS:= [\Lie(T)\times X_1(\hat
T)/\check{T}N_G(T)]\ .$$
The projection to the first factor induces a map $p:\cS\rightarrow
\cT$. The pull-back $$p^*R^*\tau:\hat \cS\rightarrow \cS$$ is given by
$$\hat \cS:=[\Lie(T)\times X_1(\hat T)/\widehat{\check{T}N_G(T)}]\ .$$
It is
trivialized (see again \cite{heinloth}, Lemma 5.6) by the $U(1)$-bundle
$$\hat \pr_2^* [U(\bar L)/\widehat{\check{T}N_G(T)}]\rightarrow [\Lie(T)\times X_1(\hat
T)/\widehat{\check{T}N_G(T)}]\ ,$$ where $$\hat \pr_2:[\Lie(T)\times X_1(\hat T)/\widehat{\check{T}N_G(T)}]\rightarrow [X_1(\hat T)/\widehat{\check{T}N_G(T)}]$$ denotes the projection.

\subsubsection{}\label{phidef}

Let $$\bar \cS:=[\Lie(T)\times X_1(\hat T)/\hat W]$$ and
$m:\cS\rightarrow \bar \cS$ be the (non-representable) map of stacks
induced by the projection $h:\check{T}N_G(T)\rightarrow \hat W$.
Finally we consider the stack $\cI:=[X_1(\hat
T)/\hat W]$ and
let $r:\bar \cS\rightarrow [ X_1(\hat
T)/\hat W]$ be induced by the projection onto the second factor.
We now consider the diagram
$$\cT\stackrel{p}{\leftarrow}\cS \stackrel{m}{\rightarrow}
\bar \cS \stackrel{r}{\rightarrow} \cI\ .$$ 
 
\subsubsection{}\label{vectb}

By our non-degeneracy assumption $\check{T}\subset \hat W$ acts freely
on $X_1(T)$. Therefore we have a diagram
$$\begin{array}{ccc}
\bar \cS&\stackrel{r}{\rightarrow}&\cI \\
\cong\downarrow&&\cong \downarrow\alpha \\
 {}[(\Lie(T)\times X_1(\hat T)/\check{T})/ W]&\stackrel{\tilde
   r}{\rightarrow}&[(X_1(\hat T)/\check{T})/W]\end{array}$$
(see \cite{heinloth}, Ex. 3.3 for the vertical isomorphisms).
Now $\tilde r$ is the projection of a  vector bundle with fibre
$\Lie(T)$.
By \ref{taud} it gives rise to a graded twist
$-\tilde \rho$ of $ 
[(X_1(\hat T)/\check{T})/W]$.
We let $\rho:=\alpha^*\tilde\rho$ be the
corresponding graded twist of $\cI$.

\subsubsection{}\label{phidef41}

The projection $r$ is now naturally $\rho$-$K$-oriented.
Moreover, $r_!:K_c( \bar
\cS)\rightarrow {}^\rho K(\cI)$ is an isomorphism.
Its inverse is the twisted Thom isomorphism.
The subscript ${}_c$ stands for proper support over $p$.
Since the fibres of $p$ are discrete this map is canonically $K$-oriented.
We define
$$\Phi:=p_!\circ (t^{-1})^*\circ m^*\circ (r_!)^{-1}:{}^\rho 
K(\cI)\rightarrow {}^{R^*\tau}K(\cT)\ .$$

\begin{prop}\label{ffds}
The map $\Phi$ is an isomorphism of groups.
\end{prop}
\proof
It suffices to show that
$$p_!\circ (t^{-1})^*\circ m^*:K_c(\bar \cS)\rightarrow
{}^{R^*\tau}K(\cT)$$
is an isomorphism.
We have an equivalence
$$\exp:\cT=[\Lie(T)/\check{T}N_G(T)]\stackrel{\sim}{\rightarrow}
[T/N_G(T)]$$
which is given by the exponential map on the level of spaces, and by
the projection $\check{T}N_G(T)\rightarrow N_G(T)$ on the level of
groups. We prove the proposition by localization over open sub-stacks of
$[T/N_G(T)]$ and the Mayer-Vietoris principle.
Since $[T/N_G(T)]$ has contractible slices we can in fact reduce to
points. 

\subsubsection{}
So let $t\in T$ and $N_G(T)_t\subset N_G(T)$ be the
stabilizer.
Then the local model is the diagram of stacks
$$[\check{T}/\check{T}N_G(T)_t]\stackrel{p}{\leftarrow}[\check{T}\times X_1(\hat
T)/\check{T}N_G(T)_t]\stackrel{m}{\rightarrow} [\check{T}\times X_1(\hat
T)/\hat W_t]\ ,$$
where $\hat W_t:=\check{T}N_G(T)_t/T\subset \hat W$, and we have identified $\exp^{-1}(\{t\})\subset \Lie(T)$ with $\check{T}$.

The restriction of the twist to the local model is given by 
$$\tau_t:[\check{T}/\widehat{\check{T}N_G(T)_t}]\rightarrow [\check{T}/\check{T}N_G(T)_t]\ .$$  It can be trivialized.
We let $N_G(T)_t$ act from the right on $\widehat{\check{T}N_G(T)_t}$ via a split
$N_G(T)_t \rightarrow \widehat{\check{T}N_G(T)_t}$.  Such a split
exists by our assumtion that the original twist $\tau$ is admissible. 
The quotient
$\widehat{\check{T}N_G(T)_t}/N_G(T)_t$ is an $U(1)$-bundle $C\rightarrow \check{T}$.
Via the left multiplication it is $\widehat{\check{T}N_G(T)_t}$-equivariant.
The bundle $C$ gives the trivialization of the restricted twist $\tau_t$.
We will use this trivialization in order to identify
$${}^{\tau_t}K([\check{T}/\check{T}N_G(T)_t])\cong
K([\check{T}/\check{T}N_G(T)_t])\ .$$

Let $U\rightarrow [X_1(\hat
T)/\widehat{\check{T}N_G(T)_t}]$ denote the restrition of $[U(\bar L)/\widehat{\check{T}N_G(T)}]$ (see \ref{gg23}) to the local model.
We form the $\check{T}N_G(T)$-equivariant $U(1)$-bundle
$V:=\hat p^*C^*\otimes \hat \pr_2^* U$, where $$\hat p:  [\check{T}\times X_1(\hat
T)/\widehat{\check{T}N_G(T)_t}]\rightarrow [\check{T}/\widehat{\check{T}N_G(T)_t}]$$ is the map induced by the first projection $p$, and $$\hat \pr_2:[\check{T}\times X_1(\hat
T)/\widehat{\check{T}N_G(T)_t}]\rightarrow [X_1(\hat T)/\widehat{\check{T}N_G(T)_t}]$$ is induced by the second projection.
 We denote by $$V^\otimes:K([\check{T}\times X_1(\hat
T)/\check{T}N_G(T)_t])\rightarrow K([\check{T}\times X_1(\hat
T)/\check{T}N_G(T)_t])$$ the operation given by the tensor product with the line bundle associated to $V$. We then must show that
$$p_!\circ V^\otimes\circ m^*:K([\check{T}\times X_1(\hat
T)/\hat W_t])\rightarrow K([\check{T}/\check{T}N_G(T)_t])$$ is an isomorphism.
The candidate for the inverse is $T_{inv}\circ (V^{\otimes})^{-1}\circ p^*$,
where $$T_{inv}:K([\check{T}\times X_1(\hat
T)/\check{T}N_G(T)_t])\rightarrow K([\check{T}\times X_1(\hat T)/\hat
W_t])$$ takes the sub-bundle of $T$-invariants. 

\subsubsection{}
In the following calculation we denote by $V$, $U$, and $C$ the complex line bundles instead of the underlying $U(1)$-bundles.
We start with the composition
$A:=T_{inv}\circ (V^{\otimes})^{-1}\circ p^*\circ p_!\circ V^\otimes\circ m^*$. All these operations can be applied on the level of vector bundles. We consider a vector bundle
 $X\rightarrow [\check{T}\times X_1(\hat T)/\hat W]$. We will show
 that $A(X)\cong X$. Indeed (where the dots indicate a
 straight-forward calculation)
\begin{eqnarray*}
A(X)&\cong&T_{inv}(\hat p^*C\otimes \hat \pr_2^* U^*\otimes p^* \circ p_!(\hat p^*C^*\otimes \hat \pr_2^* U\otimes m^*(X)))\\
&\cong&T_{inv}( \hat \pr_2^*(U^*)\otimes   \hat p^*\circ \hat p_!(\hat
\pr_2^*(U)\otimes m^*(X)))\\
&\dots&\\
&\cong&X\ .
\end{eqnarray*}


We now consider the composition 
$B:=p_!\circ V^\otimes\circ m^*\circ T_{inv}\circ (V^{\otimes})^{-1}\circ p^*$. We again calculate on the level of vector bundles
$X\rightarrow [\check{T}/\check{T}N_G(T)_t]$. We indeed have 
\begin{eqnarray*}
B(X)&\cong&p_!(\hat p^*C^*\otimes \hat \pr_2^*(U)\otimes m^*\circ
T_{inv} (\hat p^*C\otimes \hat \pr_2^* U^*\otimes p^*X))\\
&\dots&\\
&\cong& X\ .
\end{eqnarray*}
This finishes the proof of Proposition \ref{ffds}.
\hB 

\subsubsection{}

Observe (see \ref{orir})  that $R:\cT\rightarrow \cM$ is naturally $K$-oriented so that
we can consider the induction map
$$R_!:{}^{R^*\tau}K(\cT)\rightarrow {}^{\tau}K(\cM)\ .$$
\begin{prop}\label{rrf342}
On ${}^{\tau}K(\cM)$ we have
$R_!\circ R^*=\id$.
\end{prop}
\proof
We have equivalences of stacks
$$\cT\stackrel{\exp}{\cong} [T/N_G(T)]\stackrel{ind}{\cong} [G\times_{N_G(T)}T/G]$$
and $$\cM\stackrel{\hol}{\cong} [G/G]\ ,$$
where $ind$ stands for induction (see \cite{heinloth}, Ex. 3.3).
Using these equivalences we replace $R$ by the equivalent map
$$R:[(G\times_{N_G(T)}T)/G]\rightarrow [G/G]$$
(which we denote by the same symbol for simplicity).
Thus $R$ is represented by the
$G$-equivariant map
$(g,t)\mapsto gtg^{-1}$.
We can now cover $G$ by $G$-equivariant slices  $U$ on which the
twist is trivializable. We show that 
$R_!\circ R^*$ is an isomorphism, locally.
Then we argue by the  Mayer-Vietoris principle. 

\subsubsection{}

In fact we can use contractible slices. The model of the map between the slices at $t\in T$ is 
$$R_t:[(G\times_{N_G(T)_t}\Lie(T))/G]\rightarrow [(G\times_{G_t}\Lie(G_t))/G]\ .$$

We can obtain the map $R_t$  by induction of
$$S_t:[(G_t\times_{N_G(T)_t}\Lie(T))/G_t]\rightarrow [\Lie(G_t)/G_t]$$
from $G_t$ to $G$.
Thus we must show that
$(S_t)_!\circ S_t^*$ is an isomorphism.

\subsubsection{} 
By homotopy invariance and induction isomorphisms we have
$$
K([(G_t\times_{N_G(T)_t}\Lie(T))/G_t])\cong K([G_t/N_G(T)_t])\cong
K([*/N_G(T)_t])\cong R(N_G(T)_t)$$ and
$$K([\Lie(G_t)/G_t])\cong K([*/G_t]) \cong R(G_t)\ .$$
With this identification
$S_t^*:R(G_t)\rightarrow R(N_G(T)_t)$ is just the usual restriction.
In particular we know that this map is injective. In fact this is part
of the assertion of Proposition \ref{rrf342} in the untwisted
case and therefore a part of a $K$-theoretic version of the
Borel-Weyl-Bott theorem which we assume as well-known.
It therefore suffices to show that
$S_t^*(S_t)_!S^*_t=S_t^*$. Now $S_t^*(S_t)_!$ is the multiplication
with
the Euler class of the normal bundle $N$ of $S_t$.
We will show that the Euler class is equal to one.

\subsubsection{}

We consider the following diagram of $G_t$-spaces.
$$\begin{array}{ccc}
G_t\times_{N_G(T)_t}\Lie(T)&\stackrel{S_t}{\rightarrow}&\Lie(G_t)\\
\alpha\downarrow&&\beta\downarrow\\ 
G_t/N_G(T)_t&\stackrel{\delta}{\rightarrow}&*
\end{array}\ .$$
Let $N(S_t),N(\alpha),N(\beta),N(\delta)$ denote the $K_{G_t}$-classes
of normal bundles of the corresponding maps.
Then we have 
\begin{equation}\label{fuell}N(S_t)+S_t^*N(\beta)=N(\alpha)+ \alpha^*N(\delta)\ .\end{equation}
Note that $\beta,\alpha,\delta$ are submersions.
It is therefore easy to read-off the normal bundles as the inverses of
the vertical bundles.
We get
\begin{eqnarray*}
N(\alpha)&=&-\alpha^*[G_t\times_{N_G(T)_t}\Lie(T)]\\
N(\beta)&=&-\beta^*[\Lie(G_t)]\\
N(\delta)&=&-[G_t\times_{N_G(T)} (\Lie(G_t)/\Lie(T))]\ .\end{eqnarray*}
Note that
\begin{eqnarray*}
\lefteqn{\alpha^*[G_t\times_{N_G(T)_t}\Lie(T)]
+\alpha^* [G_t\times_{N_G(T)} (\Lie(G_t)/\Lie(T))]}&&\\
&=&\alpha^*[G_t\times_{N_G(T)_t}\Lie(G_t)]\\
&=&[G_t/N_G(T)_t\times \Lie(G_t)]\\
&=&S_t^*\beta^*[\Lie(G_t)]
\end{eqnarray*}
In view of \ref{fuell} this implies $N(S_t)=0$. Hence
its Euler class is the identity. This finishes the proof of Proposition \ref{rrf342} \hB

\subsubsection{}

The composition 
$$\Phi^{-1}\circ R^*:{}^{\tau}K(\cM)\rightarrow {}^\rho  K(\cI)$$ represents ${}^{\tau}K(\cM)$ as a direct summand of
${}^\rho K(\cI)$. In order to determine this summand we must
calculate  the kernel of 
$$R_!\circ \Phi: {}^\rho K(\cI)\rightarrow
{}^{\tau}K(\cM)\ .$$

\subsubsection{}
Let $\cI=\cI^{reg}\cup\cI^{sing}$ be induced by the decomposition of
$X_1(\hat T)$ into regular and singular characters.
Let $\rho_s:\hat \cI^{sing}\rightarrow \cI^{sing}$ and $\rho_r:\hat
\cI^{reg}\rightarrow \cI^{reg}$ be the corresponding restrictions of the
grading.
Then we  have a decomposition 
$${}^\rho K(\cI)\cong {}^{\rho_r} K(\cI^{reg})\oplus {}^{\rho_s}
K(\cI^{sing})\ .$$

\begin{prop}\label{p31}
We have 
$$\ker(R_!\circ \Phi)= {}^{\rho_s}
K(\cI^{sing})\ .$$
\end{prop}
\proof

\subsubsection{}\label{rch65}

The stack $\cI$ decomposes into a union of $\hat W$-orbits
$$\cI\cong \bigcup_{[\chi] \in X_1(\hat T)/\hat W}  \cI_\chi\ ,$$ where $\cI_\chi:=[[\chi]/\hat W]$.
Let $\rho_\chi$ be the restriction of $\rho$ to $\cI_\chi$.
Then we have a decomposition
$$K(\cI)\cong \bigoplus_{[\chi]\in X_1(\hat T)/\hat W}
{}^{\rho_\chi} K(\cI_\chi)\ .$$
 
We consider $\chi\in X_1(\hat T)$.
Then we must compute
$$R_!\circ \Phi_{|{}^{\rho_\chi}K(\cI_\chi) }=R_!\circ p_!\circ (t^{-1})^*\circ m^*\circ
(r_!)^{-1}_{|{}^{\rho_\chi}K(\cI_\chi)}\ .$$
In the present subsection we will show that this composition vanishes
for singular characters $\chi$. Regular characters will be discussed
later in \ref{later5}.

\subsubsection{}

For the moment, in order to set up some notation, we consider an
arbitrary character $\chi\in X_1(\hat T)$.
Let $W_\chi\subset W$ be the isomorphic image of $\hat W_\chi$ under $\hat
W\rightarrow W$. Then we have an equivalence
$$\cI_\chi\cong [(W/W_\chi)/W]\cong  [*/W_\chi]\ .$$
We let $\tilde \rho_\chi:\widehat{[*/W_\chi]}\rightarrow [*/W_\chi]$
be the induced grading. It is given by the character $\tilde \rho_\chi :
W_\chi\rightarrow \Z/2\Z$, $\tilde \rho_\chi(w):=\det_{\Lie(T)}(w)$.

\subsubsection{}\label{hchoice1}
 
In order to represent $(r_!)^{-1}_{|{}^{\rho_\chi}K(\cI_\chi)}$ we choose a $\hat W$-equivariant
section $s$ in
$$\begin{array}{ccc}
&&\Lie(T)\times [\chi]\\
&s\nearrow&\downarrow\\
{}[\chi]&=&[\chi]
\end{array} \ .$$
To this end we observe that we can choose an element $H\in \Lie(T)$ with $W_H=W_\chi$.
Then we define $s(\hat w\chi):=(\hat w H,\hat w\chi)$.

Let $s:\cI_\chi\rightarrow \bar \cS$ denote the induced map of stacks.
The normal bundle of $s$ is given by 
$[\Lie(T)\times [\chi]/\hat W]\rightarrow \cI_\chi$. Therefore $s$ is canonically
$-\rho_\chi$-$K$-oriented and we obtain a push-forward
$s_!:{}^{\rho_\chi} K(\cI_\chi)\rightarrow K(\bar \cS)$. We now
observe that
$$(r_!)^{-1}_{|{}^{\rho_\chi}K(\cI_\chi)}=s_!\ .$$

\subsubsection{}
We have the cartesian diagram
$$\begin{array}{ccc}
[[\chi]/\check{T}N_G(T)]&\stackrel{\tilde m}{\rightarrow}&\cI_\chi\\
\tilde s\downarrow&&\downarrow s\\
\cS&\stackrel{m}{\rightarrow}
&\bar \cS\end{array}\ .$$
Therefore we can rewrite
\begin{equation}\label{umschr}R_!\circ \Phi_{|{}^{\rho_\chi}K(\cI_\chi)}=R_!\circ p_!\circ \tilde s_!\circ 
(\tilde s^*t)^{-1}\circ \tilde m^*\ ,\end{equation}
The composition $R\circ p\circ \tilde s$ is represented by the 
map 
$$c:[(\hat W/\hat W_\chi)/\check{T}N_G(T)]\rightarrow \cM=[F(S^1)/G(S^1)]\ .$$ On the level of spaces it is given by 
$$\hat w \hat W_\chi\mapsto  \hat w H\ .$$
where $\hat w H$ is considered as a constant connection.
On the level of groups it is the embedding $\check{T}N_G(T)\rightarrow
G(S^1)$.
Using the equivalences $\cM\stackrel{\hol}{\cong}[G/G]$ and
$$[(\hat W/\hat W_\chi)/\check{T}N_G(T)]\cong
[(W/W_\chi)/N_G(T)]\stackrel{ind}{\cong} [(G/N_G(T)_\chi)/G]$$ we have an
equivalent representation of $c$ as $G$-equivariant 
map $$c:[(G/N_G(T)_\chi)/G]\rightarrow [G/G]$$
which is given by $[g]\mapsto \exp(gH)$. 
Note that this map factors over
$G/G_t$, so that we have a factorization of $c$ as
$$[(G/N_G(T)_\chi)/G]\stackrel{a}{\rightarrow} [(G/G_t)/G]\stackrel{b}{\rightarrow}
[G/G]\ .$$
The map $a$ is the induction of
$$a:[*/N_G(T)_\chi]\rightarrow [*/G_t]\ .$$

\subsubsection{}
From now on until the end of the proof we use the holonomy isomorphism $\cM\cong [G/G]$ in order to view $[G/G]$ as the target of $R$. Furthermore
 we will write $\tau$ for $\hol_*\tau$.
We can now write
\begin{equation}\label{jklop}R_!\circ \Phi_{|{}^{\rho_\chi}K(\cI_\chi)}= b_! \circ a_! \circ 
(\tilde s^*t)^{-1}\circ \tilde m^*\ .\end{equation}
 Note that
$b^*\tau$ is already trivial.
We choose any trivialization
$t^\prime:0\stackrel{\sim}{\rightarrow} b^*\tau$.
Then we can further write
$$R_!\circ \Phi_{|{}^{\rho_\chi}K(\cI_\chi)}= b_! \circ (t^{\prime,-1})^* \circ
a_! \circ (a^* t^\prime)^*\circ
(\tilde s^*t)^{-1}\circ \tilde m^*\ .$$

\subsubsection{} 
The embedding $[*/N_G(T)_\chi]\rightarrow [[\chi]/\check{T}N_G(T)]$
induced by $*\mapsto \chi$ and the inclusion of groups
$N_G(T)_\chi\rightarrow \check{T}N_G(T)$ is an equivalence of stacks.
 The automorphism $(\tilde s^*t)^{-1}\circ a^* t^\prime$ of the trivial twist
 can therefore be considered as a 
$N_G(T)_\chi$-equivariant line bundle on $*$. Such a line bundle 
determines a character
$\mu\in X(T)$ which is necessarily $W_\chi$-invariant.

\subsubsection{}

If $\chi$ is singular,
then $\mu$ is a character of $T$ which is singular (for $G_t$).
In this case the composition $$(a^* t^\prime)^*\circ
(\tilde s^*t)^{-1}\circ \tilde m^*:{}^{\rho_\chi}K(\cI_\chi)\rightarrow
K([*/N_G(T)_\chi])\cong R(N_G(T)_\chi)$$ produces a representations of $N_G(T)_\chi$ on
which $T$ acts by a singular characters, and
by Borel-Weyl-Bott the induction 
$a_!:K([*/N_G(T)_\chi])\rightarrow K([*/G_t])$ vanishes on
singular characters. 
This shows that $R_!\circ \Phi$ vanishes on the contributions of
singular orbits and $${}^{\rho_s}
K(\cI^{sing})\subset \ker(R_!\circ \Phi)\ .$$
The opposite inclusion $\ker(R_!\circ \Phi)\subset {}^{\rho_s}
K(\cI^{sing})$ will be shown in \ref{later5}.

\subsection{Detection of elements of twisted $K$-theory}

\subsubsection{}\label{fedef1}

We start with a twist $\tau:\hat \cM\rightarrow \cM$ which is given by
a central $U(1)$-extension $\hat G(S^1)\rightarrow G(S^1)$.
By restriction (see \ref{somp6} for the notation) we obtain a central extension
$$0\rightarrow U(1) \rightarrow \widehat{\check{T} N_G(T)}\rightarrow
\check{T}N_G(T)\rightarrow 0\ .$$ 
Since we assume that $\tau$ is admissible  this extension is trivial when restricted to
$\check{T}$ and $N_G(T)$. 

\subsubsection{}\label{bf29}

Recall (see \ref{somp6}) that $\check{T}$ acts on $X_1(\hat T)$. 
We choose a splitting homomorphism $s_0:N_G(T)\rightarrow
\widehat{\check{T} N_G(T)}$ and 
define a bilinear form
$$B:\check{T}\otimes T\rightarrow U(1)$$ by
$$B(\check{t},t):=\frac{(\check{t}\chi)(s_0(t))}{\chi(s_0(t))}\ .$$
It is easy to check that $B$ does not depend on the choice of
$\chi\in X_1(\hat T)$ and the split $s_0$. Using the form $B$ we can write the cocycle defining the
extension
\begin{equation}\label{e4r}0\rightarrow U(1) \rightarrow \widehat{\check{T} T}\rightarrow
\check{T} T\rightarrow 0\end{equation}
in the form $\omega((\check t,t),(\check t^\prime,t^\prime))= B(\check
t,t^\prime)-B(\check t^\prime,t)$.
If $s_1:\check{T}\rightarrow \widehat{\check{T}}$ is a splitting homomorphism of
$\widehat{\check{T}}\rightarrow \check{T}$, then we have
\begin{equation}\label{brel}
s_1(\check{t})s_0(t)s_1(\check{t})^{-1}=s_0(t)B(\check{t},t)\ .
\end{equation}

Note that $B$ is $W$-invariant in the sense that
$B(\check{t}^w,t^w)=B(\check{t},t)$.



\subsubsection{}\label{ctd1}

We define a subgroup $F\subset T$ by 
$$F:=\{t\in T|B(\check{t},t)=1\quad \forall
\check{t}\in\check{T}\}\subset T\ .$$
Since $B$ is $W$-invarinat the action of the Weyl group $W$ on $T$ preserves $F$.
Note that $F$
only depends on the central
extension $\hat G(S^1)\rightarrow G(S^1)$. If the twist $\tau$ is regular,
then $B$ is non-degenerated. In this case $F$ is finite.

\subsubsection{}

Let $X(\check{T})$ be the group of characters.
We have a map $b:F\backslash T\rightarrow X(\check{T})$ given by
$b(Ft)(\check{t})=B(\check{t},Ft)^{-1}$.

\begin{lem}\label{re90}
If $\tau$ is regular, then $b$ is an isomorphism.
\end{lem}
\proof
By the definition of $F$ the map $b$  is injective.
Regularity of the twist $\tau$ (see \ref{reg5412}) is equivalent to
non-degeneracy of $B$. In this case
this map $b$ is also surjective.
\hB

\subsubsection{}\label{ht13}

The restriction of  the extension (\ref{e4r}) to $\check{T}F$ is
trivialized by the choices of $s_0,s_1$ in \ref{bf29},  i.e. we have a
homomorphism $s:\check{T}F\rightarrow \widehat{\check{T}F}$ given by $s(\check{t}f):=s_1(\check{t})s_0(f)$. 

Let $\chi\in X_1(\hat T)$. Then the restriction $\chi_{|F}$
is $\check{T}$-invariant. 
Therefore we get a natural map
$s^*:X_1(\hat T)/\check{T}\rightarrow X(F)$ which is $W$-equivariant.
 \begin{lem}\label{inj33}
If $\tau$ is regular, then  the map $s^*:X_1(\hat T)/\check{T}\rightarrow X(F)$ is a bijection.
\end{lem}
\proof
Surjectivity follows from the following easy assertions.
Let $\bar \mu$ be a character of $F$. Then there exists a character of
$T$ such that $\bar \mu=\mu_{|F}$. Furthermore, using the split
$s_0:T\rightarrow \hat T$ and $\hat T\cong U(1)\times T$  we see that there exists a
character
$\chi_{\mu,s}\in X_1(\hat T)$ such that $s_0^*\chi_{\mu,s}=\mu$.
It follows that $s^*\chi_{\mu,s}=\bar \mu$.

Assume now that $s^*\chi=s^*\chi^\prime$. We write
$\chi^\prime=\chi+\lambda$ for some $\lambda\in X(T)$,
where we consider $X(T)\subset X(\hat T)$ naturally.
Then $s^*\lambda=0$. Therefore $\lambda$ pulls back from $T/F$.
In follows from Lemma \ref{re90} that there exists $\check{t}\in \check{T}$ such
that
$\lambda(t)=B(\check{t},t)$ for all $t\in T$.
Therefore, $\check{t}\chi=\chi^\prime$.
\hB

Let $X^{reg}(F)\subset X(F)$ be the subset of regular characters,
i.e. characters with trivial stabilizer in $W$. The bijection $s^*$
restricts to an identification of $W$-sets
$s^*:X^{reg}(\hat T)/\check{T}\rightarrow X^{reg}(F)$ and therefore
induces a bijection
$$X^{reg}(\hat T)/\hat W\stackrel{\sim}{\rightarrow}
X^{reg}(F)/W\ .$$

\subsubsection{}\label{thetadef}

Let us from now on assume that the twist $\tau$ is regular.
We consider the composition of maps
$$[\Lie(T)/\check{T}F]\stackrel{S}{\rightarrow} [\Lie(T)/\check{T}N_G(T)]\stackrel{R}{\rightarrow}
\cM\ .$$
The twist $$S^*R^*\tau:[\Lie(T)/\widehat{\check{T}F}]\rightarrow [\Lie(T)/\check{T}F]$$ is
trivialized by the section $s$ (see \ref{ht13}).
We let $u:0\stackrel{\sim}{\rightarrow} S^*R^*\tau$ denote the
corresponding isomorphism.
Let $q:[\Lie(T)/\check{T}F]\rightarrow [*/F]$ be induced by the projection to a
point. 
It is representable, proper, and $K$-oriented once we have fixed an
orientation of $\Lie(T)$. We can now define
a map $$\Theta:q_!\circ u^*\circ S^*\circ R^*:{}^\tau K(\cM)\rightarrow R(F)\ .$$

\subsubsection{}

Let $[\chi]\in X_1^{reg}(\hat T)/\hat W$. We trivialize the twist 
$\rho_\chi$ (see \ref{rch65})
by choosing an orientation of $\Lie(T)$ and the representative
$\chi$ of the class $[\chi]$. In fact this data orients
the projection $\Lie(T)\times \{\chi\}\rightarrow\{\chi\}$, and we 
extend this $\hat W$-equivariantly to an orientation
of $[\Lie(T)\times [\chi]/\hat W]\rightarrow \cI_\chi$.

We define the homomorphism $\sign:\hat W\rightarrow \{1,-1\}$ such
that $\sign(\hat w)=\pm 1$
depending on whether
   $\hat w:\Lie(T)\rightarrow \Lie(T)$ preserves
or reverses the orientation.

Then we consider the generator $e_{\chi}\in
{}^{\rho_\chi}K(\cI_\chi)\cong K(*)\cong \Z$ and set
$E_{\chi}:=(
R_!\circ \Phi)(e_{[\chi]})\in {}^\tau K(\cM)$. Note that
$E_{\hat w\chi}=\sign(\hat w) E_{\chi}$ for $\hat w\in \hat W$.
In the following we will often write the action of Weyl group elements as an
exponent, e.g. $\chi^w$ means $w\chi$.

\begin{theorem}\label{detec} 
We have
$$\Theta(E_{\chi})=\sum_{w\in W} \sign(w) (s^*\chi)^w \in R(F)\ .$$
\end{theorem}
\proof
By Proposition 
\ref{rrf342} we have $R^*E_{\chi}=\Phi(e_{\chi})$.
We must compute
$(q_!\circ u^*\circ S^*)(\Phi(e_{\chi}))$.
We have a diagram
$$\begin{array}{ccccc}
[\Lie(T)/\check{T}F]&\stackrel{\tilde p}{\leftarrow}&[\Lie(T)\times
X_1(\hat T)/\check{T}F]&\stackrel{\hat m}{\rightarrow}&[\Lie(T)\times
X_1(\hat T)/\hat W]\\
S\downarrow&&\tilde S\downarrow&&\|\\
\cT&\stackrel{p}{\leftarrow}&\cS&\stackrel{m}{\rightarrow}&\bar
\cS\end{array} \ ,$$
where the left square is cartesian, and the right square commutes.
We have
\begin{eqnarray*}
(q_!\circ u^*\circ S^*)(\Phi(e_{[\chi]})&=&
(q_!\circ u^*\circ S^*\circ p_!\circ (t^{-1})^*\circ m^*\circ
(r_!)^{-1})(e_{[\chi]}))\\
&=&(q_!\circ \tilde p_!\circ (\tilde p^*u)^*\circ (\tilde t^{-1})^*\circ \hat
m^*\circ (r_!)^{-1})(e_{[\chi]})\ .
\end{eqnarray*}
The trivialization $\tilde S^*t=:\tilde t:0\rightarrow \tilde p^*\circ S^*\circ R^*\tau$
is given by a $\widehat{\check{T}F}$-equivariant $U(1)$-bundle of
weight one 
$$\Lie(T)\times X_1(\hat T)\times U(1)\rightarrow \Lie(T)\times
X_1(\hat T)\ ,$$ where $F$ acts on the fibre over
$\Lie(T)\times \{\chi\}$ by $s^*\chi:=\chi\circ s:F\rightarrow U(1)$.


A small calculation shows that the  trivialization $\tilde p^* u$ is given by a
$\widehat{\check TF}$-equivariant $U(1)$-bundle of weight one $$\Lie(T)\times
X_1(\hat T)\times
U(1)\rightarrow \Lie(T)\times X_1(T)\ ,$$ where $\widehat{\check TF}$
acts on the $U(1)$-factor via its homomorphism $\widehat{\check
  TF}\rightarrow U(1)$ induced by the split $s$.

We conclude that 
\begin{equation}\label{theform3}(q_!\circ u^*\circ
  S^*)(\Phi(e_{\chi}))=\sum_{w\in W} \sign(w) (s^*\chi)^w \in R(F)\
  .\end{equation}\hB

\subsubsection{}

Let $R(F)(\sign)\subset R(F)$ be the subspace of elements satisfying
$\lambda^w=\sign(w)\lambda$ for all $w\in W$. Then the image of
$\Theta$ is contained in $R(F)(\sign)$.

\subsubsection{}\label{later5}
We now finish the proof of Proposition \ref{p31}.
After a choice of representatives $\chi\in X_1^{reg}(\hat T)$ for the equivalence classes $X_1^{reg}(\hat T)/\hat W$
the elements $e_{\chi}$, $[\chi]\in X_1^{reg}(\hat T)/\hat W$, form a
$\Z$-basis of ${}^{\rho_r}K(\cI^{reg})$.
By Theorem \ref{detec} and Lemma \ref{inj33} the composition
$\Theta\circ R_!\circ
\Phi_{|{}^{\rho_r}K(\cI^{reg})}:{}^{\rho_r}K(\cI^{reg})\rightarrow
R(F)$
is injective. This implies that
$\ker(R_!\circ \Phi)\subset {}^{\rho_s}K(\cI^{sing})$.
Thus we have finished the proof of  Proposition \ref{p31}.

\subsubsection{}\label{thebases}
By a combination of Propositions \ref{p31}, \ref{rrf342}, and
\ref{ffds} we see that
the elements
$$E_{\chi}:=R_!\circ \Phi(e_{\chi}),\quad [\chi]\in X_1^{reg}(\hat
T)/\hat W$$ form a $\Z$-basis of ${}^\tau K(\cM)$. Up to a sign this basis is
natural. The sign depends on the choice of an orientation of $\Lie(T)$. 
This finishes the proof of Theorem \ref{theo:calulation}.
\hB

\subsubsection{}

The group ${}^{\hol_*\tau }K[G/G]$ is a module over $K([*/G])\cong
R(G)$. Indeed, let $p:[G/G]\rightarrow [*/G]$ be the projection.
We use the isomorphism $R(G)\cong K([*/G])$ and the
$\cup$-product
$\cup:K([G/G])\otimes {}^{\hol_*\tau} K([G/G])\rightarrow
{}^{\hol_*\tau} K([G/G])$.
Then the module structure is given by
$U\otimes X\mapsto U\bullet X:=p^*U\cup X$, where $U\in R(G)$ and $X\in {}^{\hol_*\tau} K([G/G])$.

\subsubsection{}

Since ${}^\tau K([G/G])$ is a free $\Z$-module it embeds into its
complexification ${}^\tau K([G/G])_\C:={}^\tau K([G/G])\otimes_\Z\C$.
We let $R(G)_\C$ be the complexified group ring of $G$. The action
$\bullet$ extends to a linear action
 $\bullet:R(G)_\C\otimes_\C{}^{\tau} K([G/G])_\C\rightarrow {}^{\tau}
K([G/G])_\C$. 

\begin{theorem}\label{rrg1}
The $R(G)_\C$-module ${}^{\tau} K([G/G])_\C$ is isomorphic to a
quotient of $R(G)_\C$.
\end{theorem}
\proof
Let $\C[F]^W\cong \C[F/W]$ denote the algebra of
Weyl-invariant $\C$-valued functions on $F$. We have surjective restriction
homomorphisms
$$R(G)_\C\rightarrow R(T)_\C^W\rightarrow \C[F]^W\cong \C[F/W]\ .$$

Let $F^{reg}\subset F$ be the subset of elements with trivial stabilizer in $W$. We consider the associated vector bundle
$V:=F^{reg}\times_{W,\sign}\C\rightarrow F^{reg}/W$.
The complexification of $R(F)(\sign)$ can be identified with the
space $\Gamma(V)$ of sections of $V$ which is a $\C[F/W]$-module.
The complexification of $\Theta$ provides an injection
$$\Theta_\C:{}^{\tau} K([G/G])_\C\rightarrow
R(F)(\sign)_\C\stackrel{\sim}{\rightarrow} \Gamma(V)$$ of $\C$-vector spaces. It follows immediately from the
definition of $\Theta$, that this is a homomorphism of
$R(G)_\C$-modules. Therefore the image of $\Theta_\C$ becomes a
$\C[F/W]$-submodule. Such a submodule is completely determined by its
support $S\subset F^{reg}/W$. We define the ideal $I\subset R(G)_\C$ by the sequence
$$0\rightarrow I\rightarrow R(G)_\C\rightarrow
\C[F/W]/\C[S]\rightarrow 0\ .$$
Then have $${}^{\tau} K([G/G])_\C\cong R(G)_\C/I\ .$$
\hB




\subsection{The identity}

\subsubsection{}\label{inr41}

We consider the map of stacks $e:[*/G]\rightarrow [G/G]$ which on the
level of spaces is given by the identity element of $G$. 
Let $\tau$ be a twist of $\cM$ and $\hol_*\tau$ be the corresponding
twist of $[G/G]$. In the present subsection in order to simplify the
notation we will denote this twist simply by $\tau$.

We shall assume that $e^*\tau$ is trivial and fix a trivialization
$t:e^*\tau\stackrel{\sim}{\rightarrow}  0$. Note that we can write
$e^*\tau:[*/\hat G]\rightarrow [*/G]$ for a $U(1)$-central extension
$\hat G\rightarrow G$ of $G$. The datum of a trivialization $t$ is
equivalent to a  split $\phi:G\rightarrow\hat G$.

We assume that the twist $\sigma(G)$ (see \ref{sigmatr12}) is trivialized.
The map $e$ is $K$-oriented once we have chosen an orientation of
$\Lie(G)$ (see \ref{my1}).
Then  we define the element
$$E:=e_! t^*(1)\in {}^\tau K([G/G])\ .$$
The goal of the present subsection is an explicit calculation of $E$.
We obtain a formula for $E$ by calculating $\Theta(E)\in R(F)$, where
$\Theta$ and $F$ are as in Theorem \ref{detec}.

The element $E$ will give the identity of the  ring structure
on ${}^\tau K([G/G])$ discussed in Subsection \ref{i9i9i1}.

\subsubsection{}

We fix a positive root system $\Delta\subset X(T)$ of $(\Lie(G),\Lie(T))$ and let
$\rho:=\frac{1}{2}\sum_{\alpha\in \Delta}\alpha$.
The restriction of the split $\phi$ to the torus $T$ induces a
bijection
$\phi^*:X_1(\hat T)\rightarrow X(T)$. 
We define $\chi_{\phi,\rho}\in X_1(\hat T)$ by the condition
$\phi^*\chi_{\phi,\rho}=\rho$. According to Theorem
\ref{theo:calulation}
this character determines a basis element $E_{\chi_{\phi,\rho}}\in
{}^{\tau}K([G/G])$. 

Note that the orientation of $\Lie(G)$ induces an
orientation of $\Lie(T)$, since $\Lie(G)\cong \Lie(T)\oplus
\Lie(G)/\Lie(T)$, and the choice of the positive root system $\Delta$
fixes a complex structure (and hence an orientation) on
$\Lie(G)/\Lie(T)$. This fixes the sign of the basis element. 
\begin{theorem}\label{theunit}
We have $E=E_{\chi_{\phi,\rho}}$.
\end{theorem}

\subsubsection{}
The remainder of the present subsection is devoted to the proof of
this formula. The character $\rho$ of $T$ determines an associated  $G$-equivariant
line bundle $G\times_{T,\rho}\C\rightarrow G/T$.
Let
$x_\rho\in K([(G/T)/G])$ denote the $K$-theory element represented by
this bundle. 
Let $b:[(G/T)/G]\rightarrow [*/G]$ be induced by the projection
$G/T\rightarrow *$.
The choice of $\Delta$ induces a $G$-equivariant complex structure on
the tangential bundle
$T(G/T)$. This gives the $K$-orientation of $b$.
 The following is well-known in representation theory and a consequence
of the Borel-Weyl-Bott theorem. 
\begin{lem}
We have
$$b_!x_\rho=1\in K([*/G])\cong R(G)\ .$$
\end{lem}

\subsubsection{}

We consider now the diagram
$$\begin{array}{ccc}
[*/G]&\stackrel{e}{\rightarrow}&[G/G]\\
b\uparrow&&R\uparrow\\
{}[(G/T)/G]&\stackrel{c}{\rightarrow}&[T/N_G(T)]\end{array}\ .$$
The lower horizontal map is given by
$c:[(G/T)/G]\stackrel{\sim}{\rightarrow}[*/T]\stackrel{\tilde
  c}{\rightarrow} [T/N_G(T)]$,
where $\tilde c$ is given by the identity of $T$ on the level of spaces, and
by the embedding $T\mapsto N_G(T)$ on the level of groups.
We see that
$$e_! t^*(1)=e_!t^*(b_! x_\rho)=e_!b_!b^*(t)^*(x_\rho)=R_!c_!
b^*(t)^*(x_\rho)\ .$$ 
Under $[(G/T)/G]\stackrel{\sim}{\rightarrow}[*/T]$ the element
$x_\rho$ corresponds to $\rho\in R(T)\cong K([*/T])$.
The trivialization $t$ induces some trivialization
$v:\tilde c^*R^*\tau\stackrel{\sim}{\rightarrow} 0$.
Actually $\tilde c^*R^*\tau:[*/\hat T]\rightarrow [*/T]$, and
$v$ is induced by a split $\phi:T\rightarrow \hat T$.
At the moment we can take an arbitrary split, but 
later (see \ref{dedu4}) it will be important that $\phi$ is the restriction of
$\phi:G\rightarrow \hat G$ (see \ref{inr41}).
Using the notation of \ref{thetadef} we see that we 
must calculate $$q_!\circ u^*\circ S^* \circ \tilde c_!\circ 
v^*(\rho)\in R(F)\ .$$

\subsubsection{}
In the next paragraph \ref{trcw} we perform a longer calculation of
certain pull-back diagrams. Since we will use this result in later
subsections with different input we will state it in a general form. For the purpose of the
present subsection the  symbols have the following meaning.
\begin{enumerate}
\item $X:=\Lie(T)$, $Q:=\check{T}N_G(T)$
\item $Z:=*$, $L:=T$
\item $Y:=\Lie(T)$, $H:=\check{T}F$
\end{enumerate}
We have natural homomorphisms $L\rightarrow Q\leftarrow H$.
Furthermore, we have a central extension
$\hat Q\rightarrow Q$ given by $\widehat{\check{T}N_G(T)}\rightarrow
\check{T}N_G(T)$  such that its restrictions $\hat H\rightarrow H$ and
$\hat L\rightarrow L$ are trivialized by sections $s$ and $\phi$.
We consider the pull-back diagram
$$\begin{array}{ccccc}P&\stackrel{
  a}{\rightarrow}&[Y/H]&\stackrel{q}{\rightarrow}&[*/H]\\ d\downarrow&& S\downarrow&&\\ {}[Z/L]&\stackrel{\tilde c}{\rightarrow}&[X/Q]&&\end{array}\ .$$
Over the left upper corner $P$ we have two
trivializations of the twist $a^*S^*R^*\tau\cong d^*\tilde
c^*R^*\tau$, namely $a^*s$ and $d^*\phi$. The composition $U:=a^*u\circ
d^*v$ is an automorphism of the trivial twist and therefore a line
bundle over $P$. Using $S^*\circ \tilde c_!=a_!\circ d^*$ we see that we must calculate
$$q_! \circ a_!([U]\cup (d^*\rho))\in R(F)\ ,$$
where $[U]\in K(P)$ is the element represented by $U$.

\subsubsection{}\label{trcw}

We consider the following diagram of stacks
$$[Z/L]\rightarrow [X/Q]\leftarrow [Y/H]\ .$$
We make the pull-back $P$ explicit and get
$$\begin{array}{ccc}
[(Z\times Q)\times_X(Y\times Q)/L\times H\times Q]&\rightarrow&[Y/H]\\
\downarrow&&\downarrow\\
{}[Z/L]&\rightarrow&[X/Q]
\end{array}\ .$$
Here the   $(Z\times Q)\times_X (Y\times Q)\subset Z\times Q\times
Y\times Q$ is defined by the equation
$g_1 z=g_2 y$, where $(z,g_1,y,g_2)\in Z\times Q\times
Y\times Q$.
The action is given by
$(l,h,g)(z,g_1,y,g_2)=(lz,gg_1l^{-1},hz,gg_2h^{-1})$.
The left vertical map is given by $(z,g_1,y,g_2)\mapsto z$ on the level of
spaces, and by $(l,h,g)\mapsto l$ on the level of groups.
The upper horizontal map is similar.

We have further pull-backs
$$\begin{array}{ccc}
[Z\times Q/L\times \hat Q]&\rightarrow&[X/\hat Q]\\
\downarrow&&\downarrow\\
{}[Z/L]&\rightarrow&[X/Q]
\end{array}\ ,\quad \begin{array}{ccc}
[Y\times Q/H\times \hat Q]&\rightarrow&[X/\hat Q]\\
\downarrow&&\downarrow\\
{}[Y/H]&\rightarrow&[X/Q]\ .
\end{array}$$
We give more details for the left square.
The action is given by 
$(l,\hat g)(z,g)=(lz,\hat ggl^{-1})$.
The left vertical map is given by $(z,g)\mapsto z$ on the level of
spaces, and by $(l,\hat g)\mapsto l$ on the level of groups.
The upper horizontal map is given by
$(z,g)\mapsto gz$ on the level of spaces, and by $(l,\hat g)\mapsto
\hat g $ on the level of groups.
We pull-back further:
$$\begin{array}{ccc}
[(Z\times Q)\times_X(Y\times Q)/L\times H\times \hat Q]&\rightarrow&[Z\times Q/L\times \hat Q]\\
\downarrow&&\downarrow\\
{}[(Z\times Q)\times_X(Y\times Q)/L\times H\times Q]&\rightarrow&[Z/L]
\end{array}\ .$$
The split of $\phi:\hat L\rightarrow L$ gives rise to an action
of $L$ on $\hat Q$. We form the $U(1)$-bundle of weight one
$$[Z\times \hat Q/L\times \hat Q]\rightarrow [Z\times Q/L\times \hat
Q]\ ,$$
where the action is given by $(l,\hat g)(z,\hat g_1)=(lz,\hat g\hat
g_1\phi(l^{-1}))$.
We consider the pull-back

$$\begin{array}{ccccc}
\cL&:=&[(Z\times \hat Q)\times_X(Y\times Q)/L\times H\times \hat
Q]&\rightarrow&[Z\times \hat Q/L\times \hat Q]\\
&&\downarrow&&\downarrow\\
&&{}[(Z\times Q)\times_X(Y\times Q)/L\times H\times \hat
Q]&\rightarrow&[Z\times Q/L\times \hat Q]\\
\end{array}\ ,$$
with the action
$(l,h,\hat g)(z,\hat g_1,y,g_2)=(lz,\hat g\hat g_1 \phi(l^{-1}),\hat g g_2
h^{-1})$.
A similar construction with $[Y/H]$ gives the $U(1)$-bundle
$$\cH:=[(Z\times  Q)\times_X(Y\times \hat Q)/L\times H\times \hat
Q]\rightarrow [(Z\times Q)\times_X(Y\times Q)/L\times H\times \hat
Q]$$
with the action
$(l,h,\hat g)(z,g_1,y,\hat g_2)=(lz,\hat gg_1 l^{-1},\hat g \hat g_2
s(h^{-1}))$.
The bundle of fibrewise $U(1)$-isomorphisms $$\Hom(\cL,\cH)\rightarrow [(Z\times Q)\times_X(Y\times Q)/L\times H\times \hat
Q]$$ admits an action of
$L\times H\times Q$ in a natural way.

We simplify the description of the bundle
$$[\Hom(\cL,\cH)/L\times H\times Q]\rightarrow [(Z\times  Q)\times_X(Y\times Q)/L\times H\times 
Q]\ .$$
We first consider the special case
$$[\Hom(\hat Q\times Q,Q\times \hat Q)/Q]\rightarrow [Q\times Q/Q]$$
where $Q$ acts diagonally.
This bundle is equivalent to
$$\Hom(\hat Q,Q\times U(1))\rightarrow Q\ .$$
In our case the bundles come with an action of $L\times H$.
In the simplified picture the action is given by 
$(l,h)(\psi)(\hat g)=\psi(s(h)\hat g \phi(l^{-1}))$ over
$(l,h)(g)=hgl^{-1}$ on $Q$.
We obtain the equivalent description
in the form
$$[(Z\times \Hom(\hat Q,Q\times U(1)))\times_X Y/L\times
H]\rightarrow
[(Z\times Q)\times_XY/L\times H]\ ,$$
where
the action is given by
$(l,h)(z,\phi,y)=(lz,(l,h)(\phi),hy)$ over
$(l,h)(z,g,y)=(lz,hgl^{-1},hy)$, and the subscript $\times_X$
stands for the relation $gz=y$.

\subsubsection{}

We use the calculation \ref{trcw} in order to compute $U$.
The result is
\begin{eqnarray*}
U&\cong &[\widehat{\check{T}N_G(T)}^*\times_{\Lie(T)} \Lie(T)/T\times
\check{T}F]\\
P&\cong &[\check{T}N_G(T)\times_{\Lie(T)} \Lie(T)/T\times \check{T}F]
\end{eqnarray*}
The condition $\times_{\Lie(T)}$ stands for
$\hat n1=l$, where $(n,l)\in \widehat{\check{T}N_G(T)}^*\times
\Lie(T)$. We describe the action in the case of $U$, where it is given
by
$$(t,\check{t}f)(\hat n,l)=( s(\check{t}f)\hat n
\phi(t^{-1}),\check{t}l)\ .$$
Using the condition we can simplify the description of $U$ and $P$ to
\begin{eqnarray*}
U&\cong &[\widehat{\check{T}N_G(T)}^*/T\times
\check{T}F]\\
P&\cong &[\check{T}N_G(T)/T\times \check{T}F]
\end{eqnarray*}
If we tensorize $L$ with $d^*\rho$, then we get
$U(\rho):=[\widehat{\check{T}N_G(T)}/T\times
\check{T}F]$, where the action is now given by 
$$(t,\check{t}f)(\hat n)= s(\check{t}f)\hat n
\phi(t^{-1})\rho(t)\ .$$

\subsubsection{}

We write
$$P:=\sqcup_{w\in W} [\check{T}Tw/T\times \check{T}F]\ .$$
For $w\in W$ let   $b_w:[\check{T}Tw/T\times \check{T}F]\rightarrow
[*/F]$
be the restriction of $q\circ a$ to the corresponding component.
Then we must calculate
$(b_w)_![U_w(\rho)]\in R(F)$, where $U_w(\rho)$ is the restriction of
$U(\rho)$ to $\check{T}Tw$.
The inclusion
$[*/F]\rightarrow [\check{T}Tw/T\times \check{T}F]$ given by
$*\mapsto w$ on the level of spaces, and by $f\mapsto (f^{w^{-1}},f)$ on the
level of groups is an isomorphism of stacks. Let $\omega(f):=s(f)\phi(f^{-1})$. Then
the restriction of
$U_w(\rho)$ to $[*/F]$ is the character $\rho(f^{w^{-1}})\omega(f)$.
We obtain
\begin{equation}\label{erj9}\Theta(E)=\sum_{w\in W} \sign(w) \rho^w \omega\ ,\end{equation}
where the sign is obtained by determining the orientation of $b_w$ for
all $w\in W$ back-tracing the definitions.
\subsubsection{}\label{dedu4}

Recall that $s_{|F}:F\rightarrow \widehat{\check{T}T}$ comes as a
restriction of a split $s_0:T\rightarrow \hat T$. For the moment we
can choose $s_0:=\phi$. Then we have $\omega=0$ and therefore
$\Theta(E)=\sum_{w\in W} \sign(w)(s^*\chi_{\phi,\rho})^w$.  
This implies by Theorem \ref{detec} that
$E=E_{\chi_{\phi,\rho}}$. \hB

\subsection{The anti-diagonal}

\subsubsection{}

In the present subsection we again assume that the twist $\sigma(G)$
(see \ref{sigmatr12}) 
is trivial.
We consider the homomorphism
$$G\rightarrow G\times G\ ,\quad g\mapsto (g,g)$$
and the $G$-equivariant map
$$G\mapsto G\times G \ ,\quad g\mapsto (g,g^{-1})\ ,$$ where $G$ acts on $G$ and
$G\times G$ acts on $G\times G$ by conjugation.
In this way we arrive at a map of stacks
$$\delta:[G/G]\rightarrow [G\times G/G\times G]$$
which we call the anti-diagonal map.

\subsubsection{}

Let $\tau$ be a twist of $\cM$. We assume that $\tau$ is regular,
admissible (\ref{reg5412}), and odd  (\ref{oddef1}). The twist $\tau$ induces a twist  $\hol_*\tau$ of
$[G/G]$. In the present subsection we will simplify the notation and
write $\tau$ for this twist.  Let $\pr_i:[G\times G/G\times G]\rightarrow
[G/G]$ denote the projections. We obtain a regular and admissible twist
$\sigma:=\pr_1^*\tau+\pr_2^*\tau$ of $[G\times G/G\times G]$.
It is given by the central extension
$\widehat{(G\times G)(S^1)}\rightarrow (G\times G)(S^1)$, where
$\widehat{(G\times G)(S^1)}:=\hat G(S^1)\times \hat G(S^1)/U(1)$, and
the quotient is by the diagonal action.

\subsubsection{}

Note that $\pr_1\circ \delta=\id$ and $\pr_2\circ \delta=I$ (see \ref{oddef1}).
Since $\tau$ is odd, it follows that $\delta^*\sigma\cong
\tau+I^*\tau\cong 0$.
We choose a
trivialization $t:\delta^*\sigma\stackrel{\sim}{\rightarrow} 0$.
Note that $\delta$ is representable and proper.
We consider the composition
$$[G\times G/G\times
G]\cong[G/G]\stackrel{\delta}{\rightarrow}[G\times
G/G\times G]\rightarrow [*/G\times G]\ .$$
Here the action of $G\times G$ on $G\times G$ in the left $[G\times G/G\times
G]$ is given by $(g_1,g_2)(h_1,h_2)=(g_1h_1g_1^{-1},g_2h_2)$.
This composition  is equivalent to  the product of $[G/G]\rightarrow [*/G]$ and
$*\rightarrow [*/G]$. By \ref{my1} and \ref{my2} and the assumption
that $\sigma(G)=0$ these maps are
$K$-oriented if we choose an orientation of $\Lie(G)$. Then by \ref{twoofthree} also 
$\delta$ is $K$-oriented. 
Note that the orientation of $\delta$ is even independent of the choice of the orientation of $\Lie(G)$.

We can now define 
$$D:=\delta_! t^*(1)\in {}^\sigma K([G\times G/G\times G])\ .$$
The goal of the present subsection is an explicite calculation of $D$.

\subsubsection{}\label{om34}

The class $D$ depends on the choice of the trivialization
$t$. In the following paragraph we will define a character $\omega\in
X(T)$ which encodes the choice of $t$.

We consider the diagram
$$\begin{array}{ccc}[F(S^1)/G(S^1)]&\stackrel{\tilde \delta}{\rightarrow}&[F(S^1)\times F(S^1)/G(S^1)\times G(S^1)]\\
\cong\downarrow&&\cong \downarrow\\
{}[G/G]&\stackrel{ \delta}{\rightarrow}&[G\times G/G\times G]
\end{array}\ ,$$
where the vertical equivalences are given by the holonomy maps.
The pull-back of $\sigma$ to the right upper corner is given by 
$$[F(S^1)\times F(S^1)/\widehat{G(S^1)\times G(S^1)}]\rightarrow
[F(S^1)\times F(S^1)/G(S^1)\times G(S^1)]\ ,$$ where
$\widehat{G(S^1)\times G(S^1)}\cong \hat G(S^1)\times \hat
G(S^1)/U(1)$ (diagonal action), and $\hat G(S^1)\rightarrow G(S^1)$
defines $\tau$. 

Let $i:S^1\rightarrow S^1$ be the inversion map. The map $\tilde
\delta$ is given by $A\mapsto (A,i^*A)$ and $g\mapsto (g,i^*g)$ in the
level of spaces and groups.

The lift of the pull-back $\delta^*\sigma$ to the left upper corner is
represented by the central extension
$\widehat{G(S^1)}^d\rightarrow G(S^1)$,
which is obtained as the restriction of 
$\widehat{G(S^1)\times G(S^1)}\rightarrow G(S^1)\times G(S^1)$ via the
embedding $(\id,i^* ):G(S^1)\rightarrow G(S^1)\times G(S^1)$.
The trivialization $t$ is now given by a split $\phi:G(S^1)\rightarrow
\widehat{G(S^1)}^d$.

Let $s_0:T\rightarrow \hat T$ denote the split which
was chosen in \ref{bf29}. Then we define another split
$s_0^\prime:T\rightarrow \hat T$ by the condition that
$\phi(t)=[(s_0(t),s_0(t)^\prime)]$, where we consider $T\subset
G(S^1)$, $\hat T\subset
\hat G(S^1)$, and the bracket on the right-hand side denotes the class in
$\hat T\times_T \hat T/U(1)\subset \widehat{G(S^1)}^d$.

We define $\omega\in X(T)$ by the condition that
$s_0^\prime(t)=\omega(t)s_0(t)$, $t\in T$.

\subsubsection{}

Let $\widehat{T\times T}\rightarrow T\times T$ be the central
extension obtained by the restriction of $\widehat{(G\times G)(S^1)}$ to
$T\times T\subset (G\times G)(S^1)$, i.e
$\widehat{T\times T}=\hat T\times \hat T/U(1)$ (diagonal action).
We define a map
$$\times :X_1(\hat T)\times X_1(\hat T)\rightarrow
X_1(\widehat{T\times T})$$
by $(\chi\times \chi^\prime)([\hat t,\hat t^\prime])=\chi(\hat
t)\chi^\prime(\hat t^\prime)^{-1}$.
This map is a bijection.

\subsubsection{}
The affine Weyl group of $G\times G$ is isomorphic to $\hat W\times
\hat W$. The product $\times $ is equivariant in the sense that
$\hat w\chi\times \tilde{ \hat w}^\prime\chi^\prime=(\hat w,\hat
w^\prime)(\chi\times \chi^\prime)$, where
$\tilde{\dots} :\hat W\rightarrow \hat W$ is the automorphism
$\widetilde{\check{t}w}:=\check{t}^{-1} w$.
In particular, the product $\times$ identifies the subset
$X_1^{reg}(\hat T)\times X_1^{reg}(\hat T)$ with
$X^{reg}_1(\widehat{T\times T})$.
Furthermore we see that
it descends to a map of orbits
$$\bar \times:X_1(\hat T)/\hat W\times X_1(\hat
T)/\hat W\stackrel{\sim}{\rightarrow} X_1(\widehat{T\times T})/\hat W\times
\hat W\ .$$

\subsubsection{}

We apply Theorem \ref{theo:calulation} to the group $G\times G$. It
provides basis elements  $E_{\kappa}\in {}^{\sigma}K([G\times G/G\times
G])$ labeled by $[\kappa]\in X_1^{reg}(\widehat{T\times T})/(\hat
W\times \hat W)$. The orientation of $\Lie(T\times T)\cong
\Lie(T)\oplus\Lie(T)$ induced by a choice of an
orientation of $\Lie(T)$ is independent of this choice. This fixes the
signs of the basis elements.

\begin{theorem}\label{thedi2}
We have
$$D=\pm \sum_{[\chi]\in X_1^{reg}(\hat T)/\hat W} E_{\chi \times (\chi\omega)}\ .$$
\end{theorem}
We will explain the origin of the sign during the proof \ref{signfix}.
The idea of the proof is  to apply \ref{detec}.

\subsubsection{}

We consider the cartesian diagram
$$\begin{array}{ccc}
[T/N_G(T)]&\stackrel{j}{\rightarrow}&[T\times T/N_G(T)\times N_G(T)]\\
R\downarrow&&R\times R\downarrow \\
{}[G/G]&\stackrel{\delta}{\rightarrow}&[G\times G/G\times G]\end{array}\ ,$$
where $j$ is given by the anti diagonal $t\mapsto (t,t^{-1})$ on the level of spaces, and
by the diagonal $t\mapsto (t,t)$ on the level of groups. 
The $K$-orientation of $\delta$ induces a $K$-orientation of $j$.
We now observe that
$1=R_!R^*(1)=R_!(1)$. Therefore we have
$D=(R\times R)_!j_! R^*(t)^*(1)$.

\subsubsection{}\label{formde3}

Let $B:\check{T}\times T\rightarrow U(1)$ be the bilinear form (see \ref{bf29})
defining
$\widehat{\check{T}T}$. Then the form  
$\tilde B:(\check{T}\times \check{T})\times (T\times T)\rightarrow U(1)$
which defines $\widehat{\check{T}T\times \check{T}T}$ is given by
$$\tilde B((\check
t,\check{t}^\prime)(t,t^\prime))=B(\check{t},t)B(\check{t}^\prime,t^\prime)^{-1}\
.$$ 
Let $F\subset T$ be the subgroup introduced in \ref{ctd1}. Then we
must consider $F\times F\subset T\times T$, correspondingly.

\subsubsection{}

We now consider the diagram
\begin{equation}\label{locref}\begin{array}{ccc}
P&\stackrel{\alpha}{\rightarrow}&[T/N_G(T)]\\
p\downarrow&&j\downarrow\\
{}[T\times T/F\times F]&\stackrel{S}{\rightarrow}&[T\times T/N_G(T)\times N_G(T)]\\
q\downarrow&&\\
{}[*/F\times F]&&
\end{array} \ ,\end{equation}
where $P$ is the pull-back.
We choose a trivialization
$u:0\rightarrow S^* (R\times R)^*\sigma$.
Then the injection $\Theta$ (Theorem \ref{detec}) is given by
$$\Theta=q_!u^*S^*(R\times R)^*:{}^{\sigma} K([G\times G/G\times
G])\rightarrow R(F\times F)\ .$$
Therefore we must calculate
$$\Theta(D)=q_!u^*S^*j_!R^*(t)(1)=q_!p_! (U)\ ,$$ 
where $U\rightarrow P$ is the line bundle which corresponds to
the automorphism $\alpha^*R^*(t)\circ p^*u$ of the trivial twist.
Note that the $K$-orientation of $p$ is induced from the $K$-orientation of
$j$, and the $K$-orientation of $q$ is given by the orientation of
$\Lie(T\times T)$.

\subsubsection{}

We write the pull-back diagram (\ref{locref}) in the following
equivalent form
\begin{equation}\label{locref1}\begin{array}{ccc}
P&\stackrel{\alpha}{\rightarrow}&[\Lie(T)/\check{T}N_G(T)]\\
p\downarrow&&d\downarrow\\
{}[\Lie(T)\times \Lie(T)/\check{T}F\times \check{T}F]&\stackrel{\beta}{\rightarrow}&[\Lie(T)\times \Lie(T)/\check{T}N_G(T)\times \check{T}N_G(T)]\\
q\downarrow&&\\
{}[*/F\times F]&&
\end{array} \ .\end{equation}
The map $d$ is given by $d(l):=(l,l^{-1})$ on the level of spaces, and
by $d(m):=(m,\tilde m)$ on the level of groups.
Here $m\mapsto \tilde m$ denote the automorphism
$\check{T}N_G(T)\rightarrow \check{T}N_G(T)$ given by
$\widetilde{\check{t}u}=\check{t}^{-1}u$.

\subsubsection{}\label{se428}

The twist
$(R\times R)^*\sigma$ is given by the central extension
$$\widehat{\check{T}N_G(T)\times \check{T}N_G(T)}\rightarrow
\check{T}N_G(T)\times \check{T}N_G(T)\ .$$
By restriction via $d$ we obtain a central extension
$$\widehat{\check{T}N_G(T)}^d\rightarrow \check{T}N_G(T)\ .$$
The restriction of the split
$\phi:G(S^1)\rightarrow\widehat{G(S^1)}^d$ (see \ref{om34})  
 induces a split
$\phi:\check{T}N_G(T)\rightarrow \widehat{\check{T}N_G(T)}^d\subset \widehat{\check{T}N_G(T)\times \check{T}N_G(T)}$.

The restriction $\widehat{\check{T}F\times \check{T}F}$ of $
\widehat{\check{T}N_G(T)\times \check{T}N_G(T)}$ to $\check{T}F\times
\check{T}F$ also admits a section $\tilde s:\check{T}F\times
\check{T}F\rightarrow \widehat{\check{T}F\times \check{T}F}$.
We will write $s_1,s_2:\check{T}F\rightarrow
\widehat{\check{T}F}\subset \widehat{\check{T}F\times \check{T}F}$ for
the restrictions of $\tilde s$ to the first and the second factors.
We can and will assume that $[s_1(t),s_2(t)]=[s_0(t),s_0^\prime(t)]$,
$t\in T$.

\subsubsection{}

In order to calculate the bundle $U\rightarrow P$ explicitly we will
employ the general calculation \ref{trcw}. The symbols in \ref{trcw} now have
the following meanings. 
\begin{enumerate}
\item
$X:=\Lie(T)\times \Lie(T)$, $Q:=\check{T}N_G(T)\times \check{T}N_G(T)$
\item 
$Y:=\Lie(T)\times \Lie(T)$, $H:=\check{T}F\times \check{T}F$
\item
$Z:=\Lie(T)$, $L:=\check{T}N_G(T)$.
\end{enumerate}
We observe that the restriction
$\hat H\rightarrow H$ is trivialized by a section $\tilde s:H\rightarrow \hat
H$, and that we also have a section $\phi:L\rightarrow \hat L$.

\subsubsection{}

The calculation of \ref{trcw} gives now the following description
of the line bundle
$U\rightarrow P$.
\begin{eqnarray*}
U&\cong &[\Lie(T)\times (\widehat{\check{T}N_G(T)\times
  \check{T}N_G(T)})^*\times_{\Lie(T)\times \Lie(T)}\Lie(T)\times
\Lie(T)]/\check{T}N_G(T)\times \check{T}F\times \check{T}F]\\
P&\cong&[\Lie(T)\times \check{T}N_G(T)\times
  \check{T}N_G(T)\times_{\Lie(T)\times \Lie(T)}\Lie(T)\times
\Lie(T)]/\check{T}N_G(T)\times \check{T}F\times \check{T}F]\ .
\end{eqnarray*}
Evaluating the condition $\times_{\Lie(T)\times \Lie(T)}$
we can simplify this description
to
\begin{eqnarray*}
U&\cong &[\Lie(T)\times (\widehat{\check{T}N_G(T)\times
  \check{T}N_G(T)})^* ]/\check{T}N_G(T)\times \check{T}F\times \check{T}F]\\
P&\cong &[\Lie(T)\times \check{T}N_G(T)\times
  \check{T}N_G(T) ]/\check{T}N_G(T)\times \check{T}F\times \check{T}F]\ .
\end{eqnarray*}
Next we restrict to the section
$\Lie(T)\times \{1\}\times \check{T}N_G(T)\subset \Lie(T)\times \check{T}N_G(T)\times
  \check{T}N_G(T)$.
We get
\begin{eqnarray*}
U&\cong &[\Lie(T)\times (\widehat{\check{T}N_G(T)}^r)^* ]/\check{T}F\times \check{T}F]\\
P&\cong&[\Lie(T)\times \check{T}N_G(T) ]/\check{T}F\times \check{T}F]\ ,
\end{eqnarray*}
where $\widehat{\check{T}N_G(T)}^r\subset\widehat{\check{T}N_G(T)\times
  \check{T}N_G(T)} $ is the preimage
of $\{1\}\times \check{T}N_G(T)$. 
We describe the action for $U$:
$$(a,b)(l,\hat m)=( al,s_2(\tilde b)\hat m s_1(a) \phi(a^{-1}))\
,$$
where we consider $s_1(a)\phi(a^{-1})\in \widehat{\check{T}N_G(T)}^r$
in the natural way. 

We now write
$$\check{T}N_G(T)=\sqcup_{w\in W} \check{T}Tw\ .$$
We further restrict to the section
$$\sqcup_{w\in W} Tw\subset \sqcup_{w\in W} \check{T}Tw$$
We get
\begin{eqnarray*}
U&\cong &[\Lie(T)\times \sqcup_{w\in W}(\hat{T}^*w)^r
]/\check{T}\times F\times F]\\
P&\cong&[\Lie(T)\times \sqcup_{w\in W}Tw  ]/\check{T}\times F\times
F]\ , 
\end{eqnarray*}
where $(\hat{T}^*w)^r\subset \widehat{\check{T}N_G(T)}^r$ is the
primage of $\{1\}\times  Tw\subset \check{T}N_G(T)\times \check{T}N_G(T)$.

Let us again describe the action on $U$:
$$(a,f_1,f_2)(l,\hat tw)=(al,s_2(f_2)s_1(f_1^w)\phi((f_1^{w})^{-1})
s_2(a^w)\hat t s_1(a^w)
\phi((a^w)^{-1})w)\ ,$$
where we  write the action of Weyl group elements as exponents.
We define a character $\lambda:\check{T}F\rightarrow U(1)$
by $$\lambda(\check{t}f):= s_1(\check{t}f)s_2(\check{t}f)\phi(\check{t}f)^{-1}\ .$$
By our choice of $s_1$ and $s_2$  we have $\lambda_{|F}\equiv 1$.

Then using the relation (\ref{brel}) in the form 
$\hat t^{s_2(a^w)}=B(a^w,t)^{-1} \hat t$ (where we must take
the inverse $B(a^w,t)^{-1}$ since we consider the dual bundle
$\hat T^*\rightarrow T$) 
we can write this action in the form
$$(a,f_1,f_2)(l,\hat tw)=(al,s_2(f_2(f_1^w)^{-1} ) 
\hat t B(a^w,t)^{-1}\lambda(a^w)^{-1} w)\ .$$
We trivialize $\hat T\cong T\times U(1)$ using a split
$\kappa:T\rightarrow \hat T$. Let $\tilde \kappa:\hat T\rightarrow
U(1)$ be the associated projection.
We can now describe
$$
U\cong [\Lie(T)\times \sqcup_{w\in W}Tw\times U(1)
]/\check{T}\times F\times F]$$
with the action
$$(a,f_1,f_2)(l,tw,z)=(al,f_2 (f_1^{w})^{-1}t,\tilde \kappa(s_2(f_2(f_1^{w})^{-1}))^{-1}
\lambda(a^w)^{-1}B(a^w,t)^{-1}z)\ .$$
For $w\in W$ we consider the homomorphism
$\kappa_w:F\times F\rightarrow F\times F$ given by
$(f_1,f_2)\mapsto (f_1^{w^{-1}},f_1f_2)$.
Let $U_w\rightarrow P_w$ be the restriction of $U$ to the component
labeled by $w\in W$. 
Then we have
\begin{eqnarray*}\kappa_w^*U_w&\cong &[\Lie(T)\times (F\backslash T) \times U(1)
/\check{T}]\otimes  [U(1)/F]\\
\kappa_w^*P_w&\cong &[\Lie(T)\times F\backslash T
/\check{T}]\times [*/F]\ ,\end{eqnarray*}

where the action of $f\in F$ on $U(1)$ is trivial,
and 
the action of $\check{T}$ on $\Lie(T)\times (F\backslash T) \times U(1)
$ is here given by
\begin{equation}\label{axw2}a(l,Ft,z)=(al,Ft,\omega(a^w)^{-1}B(a^w,Ft)^{-1}z)\ .\end{equation}
Moreover we can identify $P_w\cong [T\times (F\backslash T)]\times [*/F]$.
The projection
$v_w:\kappa_w^*P_w\rightarrow [*/F\times F]$ 
is represented on the level of groups by $f\mapsto (f,1)$.

\subsubsection{}\label{signfix}

We have the associated line bundle
$$\Psi_w:=(\Lie(T)\times F\backslash T\times \C)
/\check{T}\rightarrow T\times (F\backslash T)\ ,$$
where the action is given as in (\ref{axw2}).
We must calculate
$v_![\Psi_w]\in K(*)\cong \Z$, where
$v:T\times (F\backslash T)
\rightarrow *$ is the projection.
Let $X(\check{T})$ be the group of characters.
The map $$b_w:F\backslash T\rightarrow X(\check{T})$$
which is given by $b_w(Ft)(\check{t})=B(\check{t}^w,Ft)^{-1}$ is a
bijection by Lemma \ref{re90}.
 
We see that we can identify $\Psi_w\rightarrow T\times F\backslash T$
with the Poincar{\'e} bundle $\cP\rightarrow T\times X(\check{T})$.
It is well-known that
$\tilde v_!(\cP)=1$, where  $\tilde v:T\times X(\check{T})\rightarrow
*$
and $X(\check{T})$ is oriented as the dual torus to $T$.
We therefore get
$v_!(\Psi_w)=\pm 1$, where the sign depends on whether $b_w$ preserves the
orientation or reverses it. We have $v_!(\Psi_w)=\pm \sign(w)
v_!(\Psi_1)$,
where the sign $\pm$ only depends on the twist.

\subsubsection{}

Let $z:*\rightarrow [*/F]$ be the projection. Then we have
$z_!(1)=\sum_{\chi\in \hat F}\chi\in R(F)$.
We conclude that $(v_w)_! (U_w)=\pm \sign(w)\sum_{\chi\in \hat
  F}1\otimes \chi\in R(F\times F)$.
Now note that
$$\kappa_w^*(\sum_{\chi\in \hat F} \chi^{w^{-1}}\otimes \chi^{-1})=\sum_{\chi\in \hat
  F}1\otimes \chi\ .$$
We conclude that
$$ 
q_!\circ p_! (U)=\pm \sum_{w\in W}\sign(w)\sum_{\chi\in \hat F}\chi^{w^{-1}}\otimes
\chi^{-1}\ .$$
Note that it suffices to sum over $\hat F^{reg}\subset \hat F$.

\subsubsection{}

Let us assume that
$$D=\sum_{[\chi_1],[\chi_2]\in X_1^{reg}(\hat T)/\hat W}
a_{[\chi_1],[\chi_2]} E_{\chi_1 \times \chi_2}$$
with coefficients
$a_{[\chi_1],[\chi_2]}\in \Z$ to be determined.
Then by Theorem \ref{detec} we have by equation (\ref{theform3})
$$\Theta(D)=\sum_{[\chi_1],[\chi_2]\in X_1^{reg}(\hat T)/\hat W} a_{[\chi_1],[\chi_2]} \sum_{w_1,w_2\in W}\sign(w_1)\sign(w_2)
 (\tilde s^*(\chi_1\times \chi_2))^{(w_1,w_2)}\ .$$
We have (see \ref{se428} for notation) 
$$(\tilde s^*(\chi_1\times \chi_2))^{(w_1,w_2)}=(s_0^*\chi_1)^{w_1}\otimes 
((s_0^\prime)^*\chi_2^{-1})^{w_2}\ .$$
Now observe that $(s_0^\prime)^*\chi^{-1}=(s_0^*\chi^{-1})
\omega^{-1}=
s_0^* (\omega\chi)^{-1}$ for
$\chi\in X_1(\hat T)$.



\subsubsection{}\label{corre45}

We further rewrite 
\begin{eqnarray*}q_!\circ p_!(U)
&=&\pm \sum_{w_1,w_2\in W}\sign(w_1)\sum_{[\chi]\in 
  \hat F^{reg}/W}  (\chi^{w_2})^{w_1^{-1}}\otimes
(\chi^{w_2})^{-1}\\&=&\pm \sum_{w_1,w_2\in
  W}\sign(w_1)\sign(w_2)\sum_{[\chi]\in \hat F^{reg} /W} \chi^{w_1}\otimes
(\chi^{w_2})^{-1}\ .\end{eqnarray*}
We now want to evaluate the equality 
$\Theta(D)=q_!p_!(U)$ which explicitly has the form
\begin{eqnarray*}\lefteqn{\pm \sum_{w_1,w_2\in
  W}\sign(w_1)\sign(w_2)\sum_{[\chi]\in \hat F^{reg} /W} \chi^{w_1}\otimes
(\chi^{w_2})^{-1}}&&\\&=&\sum_{[\chi_1],[\chi_2]\in X_1^{reg}(\hat T)/\hat W} a_{[\chi_1],[\chi_2]} \sum_{w_1,w_2\in W}\sign(w_1)\sign(w_2)
(s^*_0\chi_1)^{w_1}\otimes 
(s^*_0(\omega \chi_2)^{-1})^{w_2}\ .\end{eqnarray*}
In view of the discussion in \ref{ht13} we see that
$$\pm a_{[\chi_1],[\chi_2]}=\delta_{[\chi_1],[\omega\chi_2]}\ .$$ 
This finishes the proof of Theorem \ref{thedi2}. \hB

\section{Moduli spaces and trivializations of twists}

\subsection{Stacks associated to polarized Hilbert spaces}\label{pol91}
 
\subsubsection{}\label{kk4}

Let $H=H_+\oplus H_+^\perp$ be a polarized Hilbert space, and let
$P_+$ denote the orthogonal projection onto $H_+$. In this situation 
we define the restricted unitary group
$$U_{res}(H,H_+):=\{U\in U(H)\:|\:[U,P_+]\:\mbox{trace class}\}\ .$$
We obtain the basic central extension of the connected component of
the identity of $U_{res}(H,H_+)$ in the following canonical way
(see  \cite{pressleysegal}). We
first consider the subgroup
$$E:=\{(U,Q)\in U_{res}(H,H_+)\times U(H_+)\:|\:P_+UP_+-Q\:\mbox{trace
  class}\}\ .$$
This group sits in an extension
$$0\rightarrow T\rightarrow E\rightarrow U^0_{res}(H,H_+)\rightarrow 0\ ,$$
where $T= \{Q\in U(H_+)|1-Q\:\mbox{trace class}\}$. Let $\det:T\rightarrow U(1)$ be the Fredholm
determinant. Then we define
$$\hat U^0_{res}(H,H_+):=E\times_{T,\det}U(1)\ .$$
\begin{ddd}
This is the basic central extension:
\begin{equation}\label{kk44}0\rightarrow U(1)\rightarrow \hat U^0_{res}(H,H_+)\rightarrow
U^0_{res}(H,H_+)\rightarrow 0\ .\end{equation}
\end{ddd}

\subsubsection{}\label{canspil}

Let $U^{H_+}\subset U^0_{res}(H,H_+)$ be the subgroup of isometries
fixing $H_+$. Then we have a natural split
$$
\begin{array}{ccc}
&&E\\
&\tilde s\nearrow&\downarrow\\
U^{H_+}&\rightarrow&U^0_{res}(H,H_+)\end{array}
\ .$$
given by $\tilde s(u):=(u,P_+uP_+)$.
It induces a split
\begin{equation}\label{rrcas}\begin{array}{ccc}
&& \hat U^0_{res}(H,H_+)\\
&s\nearrow&\downarrow\\
U^{H_+}&\rightarrow&U^0_{res}(H,H_+)\end{array}
\ .\end{equation}

\subsubsection{}\label{ccyy8}

For completeness we explain how one can extend the central extension
from $U^0_{res}(H,H_+)$ to $U_{res}(H,H_+)$. 
Recall that $U_{res}(H,H_+)$ comes as a (non-canonically split) extension
(semi-direct
product)
$$0\rightarrow U^0_{res}(H,H_+)\rightarrow U_{res}(H,H_+)\rightarrow
\Z\rightarrow 0\ .$$
If we choose a unitary $\sigma\in U_{res}(H,H_+)$ such that $\sigma(H_+)\subset
H_+$ and $\codim_{H_+}\sigma(H_+)=1$, then we can define
a split $\Z\ni 1\mapsto \sigma\in U_{res}(H,H_+)$.
We extend $\sigma\in \Aut(U^0_{res}(H,H_+))$ to an automorphism
$\tilde \sigma \in \Aut(\hat U^0_{res}(H,H_+))$ by
$\tilde \sigma[(U,Q),z]=[(\sigma U\sigma^{-1},Q_\sigma),z]$,
where $[U,Q]\in E$ and 
$$Q_\sigma:=\left\{\begin{array}{cc}\sigma Q\sigma^{-1}& \mbox{on
      $\sigma(H_+)$}\\ 1&\mbox{on $H_+\ominus \sigma(H_+)$}
\end{array}\right. \ .$$
We then define
$$\hat U_{res}(H,H_+):=\hat U_{res}^0(H,H_+)\rtimes \Z\ .$$
This definition depends up to isomorphism on the choice of $\sigma$ and is therefore less
canonical than the construction of the extension of the connected
component.


\subsubsection{}\label{ctw21}

The upshot of the preceeding discussion is that a polarized Hilbert
space $H=H_+\oplus H_+^\perp$ gives rise to a twist
$$\cT:[*/\hat U_{res}(H,H_+)]\rightarrow [*/ U_{res}(H,H_+)]$$
whose restriction to the identity component of the restricted unitary
group is canonical.

\subsubsection{}

We now consider the restricted 
Grassmannian 
$$Gr_{res}(H,H_+):=\{P \:\mbox{orthogonal projection}\:|\:P-P_+\:\mbox{trace class}\}\ .$$
It is a homogeneous space of $U_{res}(H,H_+)$. We have the pull-back
$$\begin{array}{ccc}
[Gr_{res}(H,H_+)/\hat
U_{res}(H,H_+)]&\stackrel{p^*\cT}{\rightarrow}&[Gr_{res}(H,H_+)/
U_{res}(H,H_+)]\\
p\downarrow&&\downarrow\\
{}[*/\hat
U_{res}(H,H_+)]&\stackrel{\cT}{\rightarrow}&[*/
U_{res}(H,H_+)]\end{array}\ .$$

\begin{lem}\label{vvcc1}
The twist  $p^*\cT$ is trivialized.
\end{lem} 
\proof
Note that $Gr_{res}(H,H_+)$ carries a determinant bundle $L\rightarrow
Gr_{res}(H,H_+)$. It was shown in \cite{pressleysegal}, 7.7.3, that the central
extension $\hat U_{res}(H,H_+)$ acts canonically\footnote{Note that
  the definition of $\hat U_{res}(H,H_+)$ and the definition of the
  action depend on the same choice of $\sigma$ (see \ref{ccyy8})} on
$L$ lifting the action of $U_{res}(H,H_+)$ on $Gr_{res}(H,H_+)$.
This line bundle gives the isomorphism
$l:0\stackrel{\sim}{\rightarrow} p^*\cT$. \hB

\subsubsection{}

Later we need the following fact. Recall that  $U^{H_+}\subset U^0_{res}(H,H_+)$
is the subgroup of transformations which fix $P_+$.
Via the canonical split (\ref{rrcas}) its acts on the fibre of $L$ over $P_+$.
\begin{lem}\label{nn13}
The group $U^{H_+}$ acts trivially on the fibre of $L$ over $P_+$.
\end{lem}

\subsection{Moduli spaces}

\subsubsection{}

Let $C$ be a non-empty  oriented one-dimensional closed Riemannian manifold.
Observe that $C$ admits a natural action of $S^1$. If $V$ is a finite-dimensional complex
Hilbert space, then we consider the Hilbert space
$H:=L^2(C,V)$. The group $S^1$ acts on $V$ in a natural way. 
We obtain a polarization $H=H_+\oplus H_+^\perp$ be taking
for $H_+$ the subspace of non-negative Fourier modes.

\subsubsection{}\label{mcf5}

Let now $G$ be a compact Lie group. Then we consider the trivial
$G$-principal bundle $P(C)$ over $C$. Furthermore, we let $F(C)$ and
$G(C)$ denote the space of flat connections and the gauge group of
$P(C)$. The group $G(C)\cong C^\infty(C,G)$ acts on $F(C)$. In this
way we obtain the stack $$\cM(C):=[F(C)/G(C)]\ .$$
This generalizes the construction given in \ref{str43}.

\subsubsection{}

Let us now assume that $G$ acts unitarily on the Hilbert space $V$.
We obtain an induced homomorphism $G(C)\rightarrow U(H)$. It is
well-known (see \cite{pressleysegal}, Sec. 6.3) that this homomorphism factors over the restricted unitary
group $U_{res}(H,H_+)$. If $\pi_1(G)$ is finite, then $G(C)$ maps to $U_{res}^0(H,H_+)$.

\subsubsection{}\label{uf49}

We have an induced map of stacks
$v:\cM(C)\rightarrow [*/U_{res}(H,H_+)]$ which can be used in order to
define the twist (see \ref{ctw21} for the definition of $\cT$)
$$v^*\cT:\hat \cM(C)\rightarrow \cM(C)\ .$$
The stack $\hat \cM(C)$ is isomorphic to
$[F(C)/\hat G(C)]$, where
$\hat G(C)\rightarrow G(C)$ is the $U(1)$-central extension obtained
as pull-back of $\hat U_{res}(H_+,H)\rightarrow U_{res}(H,H_+)$ via
the homomorphism  $G(C)\rightarrow U_{res}(H,H_+)$. If $\pi_1(G)$ is
finite, then $v^*\cT$ is canonical. In general at least the isomorphism class
of $v^*\cT$ is well-defined.

\subsubsection{}
Recall the notion of admissibility   \ref{reg5412}.
\begin{lem}
The twist $v^*\cT$ is admissible.
\end{lem}
\proof
It suffices to consider the case $C=S^1$.
Note that $G\subset G(S^1)$ preserves $H_+$. We therefore have a factorization
$G\rightarrow U^{H^+}\rightarrow U_{res}(H,H_+)$. Using \ref{canspil} we obtain a split $G\rightarrow \hat G$. In particular, then central extension
$\widehat{N_G(T)}\rightarrow N_G(T)$ is trivial.

By an explicit calculation using the definitions in the case of $SU(n)$ 
and pulling back the result to $G$ we show that the extension $\widehat{\check{T}}\rightarrow \check{T}$ is also trivial. \hB

\subsubsection{}

The extension  $\check{T}T\subset G(S^1)$ determines a bilinear form $B:\check{T}\otimes T\rightarrow U(1)$ (see \ref{bf29}). In the present subsection we calculate this form. Note that $B$ is completely determined by its derivative $b:\check{T}\times \Lie(T)\rightarrow \Lie(S^1)\cong \R$ with respect to the second entry. We consider this as a homomorphism $b:\check{T}\rightarrow \Lie(T)^*$. Note that $\check{T}$ is a lattice in $\Lie(T)$. Thus $b$ has a unique extension to a linear map $b:\Lie(T) \rightarrow \Lie(T)^*$. This is a bilinear form on $\Lie(T)$.

By an explicit calculation using the definitions we obtain the following formula.
\begin{lem}\
For $G=SU(n)$ and its  standard representation on $\C^n$ the form $b$ is given by  
$b(X,Y)=\Tr(XY)$.
\end{lem}

From this we immediately obtain the general case if $\pi_1(G)$ is finite. 
In fact, in this case $\rho:G\rightarrow SU(V)$. Let $\rho:\Lie(G)\rightarrow \Lie(V)$ be the derived representation of Lie algebras. Then we have
$$b(X,Y)=\Tr_V(\rho(X)\rho(Y))\ .$$

\subsubsection{}

The form $(X,Y)\mapsto \Tr(XY)$ on $\Lie(SU(n))$ is negative definite.
It follows that its restriction to the Lie algebra of the maximal torus of $SU(n)$ is non-degenerated. We conclude:
\begin{kor}
If $\rho:\Lie(G)\rightarrow \Lie(SU(n))$ is injective, then 
the form $b$ on $\Lie(T)$ is non-degenerated. In particular, the twist
$v^*\cT$ is regular (see \ref{reg5412}).
\end{kor}

\subsubsection{}

Let $-C$ be $C$ equipped with the opposite orientation.
Note that $\cM(-C)=\cM(C)$.
Let $\bar V$ be the complex conjugated representation to $V$.
In the following we indicate the dependence of the Hilbert spaces $H$
on $V$ by writing $H(V)$. Since $H_+$ also depends on the orientation
of $C$ we will write $H_+(V,C)$.
As in \ref{uf49} we have maps $v:\cM(C)\rightarrow [*/U_{res}(H(V),H_+(V,C))]$ and
$\bar v:\cM(C)\rightarrow [*/U_{res}(H(\bar V),H_+(\bar V,-C))]$
which induce twists $\tau$ an $\bar \tau$ of $\cM(C)$.

Let us assume that $\pi_1(G)$ is finite. 
\begin{lem}\label{canis2} 
There exists  a canonical isomorphism
$\bar \tau=-\tau$.
\end{lem}
\proof
The conjugate linear isomorphism $V\stackrel{conj}{\cong} \bar V$ induces a conjugated linear
isomorphism $H(V)\stackrel{conj}{\cong}
 H(\bar V)$ which identifies $H_+(V,C)$ with $H_+(\bar V,-C)$.
We obtain a corresponding diagram of groups
$$\begin{array}{ccc}
\hat U_{res}^0(H(V),H_+(V,C))^*&\cong&\hat U_{res}^0(H(\bar
V),H_+(\bar V,-C))\\
\downarrow&&\downarrow\\
U_{res}^0(H(V),H_+(V,C))&\cong&U_{res}^0(H(\bar V),H_+(\bar
V,-C))\end{array}\ .$$
Here for an $U(1)$-central extension $\hat A\rightarrow A$ we denote by
$\hat A^*\rightarrow A$ the opposite extension. If $A$ acts on a space
$X$, then the twist $[X/\hat A^*]\rightarrow [X/A]$ is the negative of
$[X/\hat A]\rightarrow [X/A]$. The assertion now follows. \hB 

\subsubsection{}

Let us assume a decomposition $C=C_1\cup C_2$.
This induces decompositions $H=H_1\oplus H_2$ and $H_+=H_{1,+}\oplus
H_{2,+}$. As in \ref{uf49} we consider the maps $v_i:\cM(C_i)\rightarrow
[*/U_{res}(H_i,H_{+,i})]$ and define the twists $\tau_i:= v_i^*\cT_i$.
Let $\tau:=v^*\cT$.

Let us again assume that $\pi_1(G)$ is finite.
\begin{lem}\label{sum41}
We have a canonical isomorphism
$$\pr_1^*\tau_1+\pr_2^*\tau_2=\tau\ .$$
\end{lem}
\proof
We have a natural embedding
$U_{res}(H_1,H_{1,+})\times U_{res}(H_2,H_{2,+})\rightarrow
U_{res}(H,H_+)$. Let $\widehat{U^0_{res}(H_1,H_{1,+})\times
  U^0_{res}(H_2,H_{2,+})}\rightarrow U^0_{res}(H_1,H_{1,+})\times
U^0_{res}(H_2,H_{2,+})$ be the induced central extension.
It follows from the construction \ref{kk4}
that we have a canonical identification
of $(\hat U^0_{res}(H_1,H_{+,1})\times \hat
U^0_{res}(H_2,H_{+,2}))/U(1)\cong \widehat{U^0_{res}(H_1,H_{1,+})\times
  U^0_{res}(H_2,H_{2,+})}$.

Now in general, let $A$ and $B$ be groups acting on spaces $X$ and $Y$, respectively. Furthermore let $\hat
A\rightarrow A$ and $\hat B\rightarrow B$ be  $U(1)$-central
extensions inducing twists $\alpha:[X/\hat A]\rightarrow [X/A]$ and
$\beta:[Y/\hat B]\rightarrow [Y/B]$.
Then the twist $\pr_A^*\alpha+\pr_B^*\beta$ is represented by
$[X\times Y/ ((\hat A\times \hat B)/U(1))]\rightarrow [X\times
Y/A\times B]$, where $\pr_A,\pr_B$ are the obvious projections. 

This implies the result since we have a factorization of $v$ as
$$\cM(C)\cong \cM(C_1)\times \cM(C_2)\stackrel{v_1\times
  v_2}{\rightarrow}
[*/ U^0_{res}(H_1,H_{+,1})]\times [*/U^0_{res}(H_2,H_{+,2})]\rightarrow
[*/U^0_{res}(H,H_{+})] .$$
\hB

\subsubsection{}

Let now $\Sigma$ be a two-dimensional oriented Riemannian manifold 
with non-empty boundary $\partial \Sigma$. 
 Then we consider the trivial
$G$-principal bundle $P(\Sigma)$ over $\Sigma$. Furthermore, we let $F(\Sigma)$ and
$G(\Sigma)$ denote the space of flat connections and the gauge group of
$P(\Sigma)$. The group $G(\Sigma)$ acts on $F(\Sigma)$. In this
way we obtain the stack $$\cM(\Sigma):=[F(\Sigma)/G(\Sigma)]\ .$$

\subsubsection{}\label{ev123}

Evaluation at $\partial \Sigma$ defines a homomorphism
$G(\Sigma)\rightarrow G(\partial \Sigma)$ and an equivariant map
$F(\Sigma)\rightarrow F(\partial \Sigma)$. In this way we get a map of
stacks
$$q:\cM(\Sigma)\rightarrow \cM(\partial \Sigma)\ .$$

\subsubsection{}

We fix a unitary representation $V$ of $G$. Note that $\partial
\Sigma$ is compact, oriented and Riemannian. Therefore we have a twist
$v^*\cT:\hat \cM(\partial \Sigma)\rightarrow \cM(\partial \Sigma)$.

\begin{prop}\label{mmc1}
The pull-back of twists $q^*v^*\cT$ is trivialized.
\end{prop}
\proof
The Riemannian metric 
 together
with the orientation of $\Sigma$ gives
a complex structure on $\Sigma$. A connection $A\in F(\Sigma)$ induces a
holomorphic structure $\bar \partial_A$ on the associated bundle
$V(\Sigma):=P(\Sigma)\times_{G}V$. We let
$H(A)\subset H=L^2(\partial \Sigma,V)$
denote the closure of the space of boundary values of continuous
$\bar \partial_A$-holomorphic sections of $V(\Sigma)$. Let $P(A)$
be the projection onto $H(A)$. It turns out that
$P(A)\in Gr_{res}(H,H_+)$. 
We thus
obtain a map
$$P:F(\Sigma)\rightarrow Gr_{res}(H ,H_+)\ .$$
We now observe that this map is $G(\Sigma)$-equivariant, where
$G(\Sigma)$ acts on the right-hand side via its homomorphism
$$G(\Sigma)\rightarrow G(\partial \Sigma)\rightarrow U_{res}(H ,H_+)\ .$$
Eventually we obtain the diagram of maps of stacks
\begin{equation}\label{pde431}\begin{array}{ccc}
\cM(\Sigma)&\stackrel{P}{\rightarrow}&[Gr_{res}(H
,H_+)/U_{res}(H ,H_+)]\\
q\downarrow&&p\downarrow\\
\cM(\partial \Sigma)&\stackrel{v}{\rightarrow}&[*/U_{res}(H,H_+)]
\end{array}\ .\end{equation}
The required trivialization is now given by
$$P^*l:0\stackrel{\sim}{\rightarrow} P^*p^*\cT\cong q^*v^*\cT$$
with $l$ obtained in Lemma \ref{vvcc1}.
\hB 

\subsubsection{}\label{nn2}\label{su34}

Let $C$ be a compact oriented one-dimensional Riemannian manifold. 
 We consider two orientation and metric  preserving  embeddings
$f_0,f_1:(-1,1)\times C\rightarrow \Sigma$ with disjoint images.
Then we can cut $\Sigma$ at the images $f_i(\{0\}\times C)$ and glue
again interchanging the copies. In this way we obtain a compact
oriented Riemannian two-manifold $\tilde \Sigma$ again with two embeddings
$\tilde f_0,\tilde f_1:(-1,1)\times C$. Note that there is a canonical
identification $\partial \Sigma \cong \partial \tilde \Sigma$.



\subsubsection{}

We let $F(\Sigma,\sim)\subset F(\Sigma)$ be the space of flat connections $A$ on $\Sigma$ 
with the property that $f_0^*A=f_1^*A$. We define $F(\tilde
\Sigma,\sim)\subset F(\tilde \Sigma)$ 
in a similar manner. Then we have a canonical identification
$F(\Sigma,\sim)\cong F(\tilde \Sigma,\sim)$.
 
We further define $G(\Sigma,\sim)\subset G(\Sigma)$ as the subgroup
of gauge transformations  $g$ satisfying $f_0^*g=f_1^*g$.
We define $G(\tilde \Sigma,\sim)\subset G(\tilde \Sigma)$ in a similar
manner and observe that we have a canonical
identification $G(\Sigma,\sim)\cong G(\tilde \Sigma,\sim)$.

\subsubsection{}

We get a diagram of maps of stacks
$$\begin{array}{ccccc}
&&\cM(\Sigma)&&\\
&i\nearrow&&q_\Sigma\searrow\\
{}[F(\Sigma,\sim)/G(\Sigma,\sim)]&&&&\cM(\partial \Sigma)\\
&\tilde i\searrow&&q_{\tilde \Sigma}\nearrow\\
&& \cM(\tilde \Sigma) &&
\end{array} \ .$$



\subsubsection{}

In Proposition \ref{mmc1} we have constructed  trivializations
\begin{eqnarray*}
t(\Sigma)&:&q_\Sigma^* \cT\stackrel{\sim}{\rightarrow}
0\\
t(\tilde
  \Sigma)&:&q_{\tilde \Sigma}^*\cT\stackrel{\sim}{\rightarrow} 0\ .
\end{eqnarray*}
Note that
$$i^* q_\Sigma^* \cT \cong \tilde i^*q_{\tilde \Sigma}^*\cT$$ canonically.
 
\begin{prop}\label{suin}
There exists an isomorphism of trivializations
$$i^*t(\Sigma)\cong \tilde i^*t(\tilde \Sigma)$$
of $i^*q_\Sigma^* \cT\cong
\tilde i^*q_{\tilde\Sigma}^* \cT$.
\end{prop}
\proof
Recall that the trivializations $t(h^{\Sigma})$ and $t(h^{\tilde \Sigma})$ were induced by
the equivariant bundles $P^*L$ and
$\tilde P^*L$ (see \ref{vvcc1} and \ref{pde431} for the
notation, and $\tilde{\dots}$ indicates objects associated  to $\tilde \Sigma$). It suffices to show
that $i^*P^*L$ and $\tilde i^*\tilde P^*L$ are isomorphic
as $\hat G(\Sigma,\sim)$-equivariant bundles,
where the central extension  $\hat G(\Sigma,\sim)\rightarrow
G(\Sigma,\sim)$
is defined as the restriction of $\hat G(\Sigma)\rightarrow
G(\Sigma)$, and the latter is pulled back from $\hat G(\partial
\Sigma)\rightarrow G(\partial \Sigma)$ (see \ref{uf49}). 
Therefore the following Lemma implies the proposition.

\begin{lem}
The maps
$i\circ P$ and $\tilde i\circ \tilde P$
are $G(\Sigma,\sim)$-equivariantly homotopic.
\end{lem}
\proof
Let $A\in F(\Sigma,\sim)$ and $\bar \partial_A$ and $\widetilde{\bar
  \partial_A}$ be the corresponding holomorphic structures
on $V(\Sigma)$ and $V(\tilde \Sigma)$.
It is a
by now standard trick (see \cite{bunke}) to identify the spaces
$C(\Sigma,V(\Sigma))$ and $C(\tilde \Sigma,V(\tilde \Sigma))$ in a
$G(\Sigma,\sim)$-equivariant way so that
$\delta:=\bar \partial_A-\widetilde{\bar
  \partial_A}$ is a zero-order (non-local) operator.
For $t\in [0,1]$ we can form the projection $P_t(A)$ 
onto the boundary values of solutions of
$\bar \partial_A-t\delta$. It provides the homotopy from
$P(i(A))$ to $\tilde P(\tilde i(A))$. 
\hB

\subsubsection{}

In the following paragraphs  we interpret the trivialization constructed in
\ref{mmc1} and the surgery invariance \ref{suin} in a slightly
different way. Let $\Sigma$ be a an oriented compact surface with
Riemannian metric and nonempty boundary.  
We assume a decomposition of the boundary into an ingoing and an
outgoing part: $$\partial \Sigma=\partial_i\Sigma\cup \partial_a
\Sigma\ .$$ We will equip the ingoing boundary with the orientation
which is opposite to the induced orientation.

\subsubsection{}\label{ia27}

We assume that $\pi_1(G)$ is finite and
 fix an unitary representation $V$ of $G$. 
The construction \ref{uf49} gives rise to twists
$\tau_i:\hat \cM(\partial_i \Sigma)\rightarrow \cM(\partial_i\Sigma)$
and
$\tau_a:\hat\cM(\partial_a\Sigma)\rightarrow \cM(\partial_a\Sigma)$.
We consider the correspondence
\begin{equation}\label{cne3}\cM(\partial_i
\Sigma)\stackrel{q_i}{\leftarrow}\cM(\Sigma)\stackrel{q_a}{\rightarrow}
\cM(\partial_a \Sigma)\ .\end{equation}

\begin{lem}\label{ris4}
We have a canonical isomorphism
$r:q_a^*\tau_a\stackrel{\sim}{\rightarrow} q_i^*\tau_i$.
\end{lem}
\proof
We can write $q=(q_i,q_a):\cM(\Sigma)\rightarrow
\cM(\partial_i\Sigma)\times  \cM(\partial_a\Sigma)\cong \cM(\partial
\Sigma)$. Let $\pr_i:  \cM(\partial\Sigma)\rightarrow
\cM(\partial_i\Sigma)$ and
 $\pr_a: \cM(\partial\Sigma)\rightarrow \cM(\partial_a\Sigma)$ be the
 projections. 
By Lemma \ref{sum41} and Lemma \ref{canis2} we have a canonical isomorphism
$\tau=\pr_a^*\tau_a-\pr_i^*\tau_i$, where $\tau$ is the twist of
$\cM(\partial \Sigma)$  given by \ref{uf49}. By Proposition
\ref{mmc1} we have a canonical trivialization
$q^*\tau\stackrel{\sim}{\rightarrow}0$. If we add the identity
$q_i^*\tau_i\cong q_i^*\tau_i$, then we obtain an isomorphism
$r:q_a^*\tau_a\stackrel{\sim}{\rightarrow} q_i^*\tau_i$. \hB

\subsubsection{}\label{nn11}
In this subsection we reinterpret the surgery invariance \ref{suin}.
We keep the assumption that $\pi_1(G)$ is finite.
 We consider two compact oriented surfaces $\Sigma_n$, $n=0,1$ with
non-empty boundary  which are equipped with Riemannian metrics. 
We assume product structures near the boundaries.
We assume an orientation reversing isometry
$\psi:\partial_i\Sigma_1\stackrel{\sim}{\rightarrow} \partial_a\Sigma_0$.
Then we can form the compact oriented surface
$\Sigma:=\Sigma_{0}\sharp_{\partial_a\Sigma_0\cong\partial_i\Sigma_1}\Sigma_1$ with
boundary
$\partial_i \Sigma\cong \partial_i\Sigma_0$ and 
$\partial_a \Sigma\cong \partial_a \Sigma_1$. It comes equipped  with
an induced Riemannian metric.

We extend the notation introduced in \ref{ia27} by an index
$\alpha\in\{0,1\}$ in order to indicate the surface to which the objects
belong.
We consider the following diagram
\begin{equation}\label{se4318}\begin{array}{ccccc}
\cM(\Sigma)&\stackrel{j}{\rightarrow}&\cM(\Sigma_1)&\stackrel{q_{1,a}}{\rightarrow}&\cM(\partial_a\Sigma_1)\\
i\downarrow&&\psi\circ q_{1,i}\downarrow&&\\
\cM(\Sigma_0)&\stackrel{q_{0,a}}{\rightarrow}&\cM(\partial_a\Sigma_0)&&\\
q_{0,i}\downarrow&&&&\\
\cM(\partial_i\Sigma_0)&&&&
\end{array}\ ,
 \end{equation}
where $i$ and $j$ are the canonical restriction maps.
The following is just a rewriting of \ref{suin}
\begin{kor}\label{rrc17}
We have a commutative diagram
$$\begin{array}{ccccccc}
j^*q_{1,a}^* \tau_{1,a}&\stackrel{j^*(r_1)}{\rightarrow} &j^*q_{1,i}
\tau_{1,i}&\cong&
i^*q_{0,a}^*\tau_{0,a}&\stackrel{i^*(r_0)}{\rightarrow}&i^*q_{0,i}^*\tau_{0,i}\\
\cong \|&&&&&&\cong \|\\
q_a^*\tau_a&&&\stackrel{r}{\longrightarrow}&&&q_i^*\tau_i\end{array}\
,$$
where all maps denoted by $\cong$ are canonical identifications.
\end{kor}

\subsection{The product}\label{i9i9i1}

\subsubsection{}

In the present subsection we shall assume that $\pi_1(G)$ is
finite. This implies that the twists constructed in \ref{uf49}
and their trivializations \ref{mmc1} are canonical.

Moreover we assume that the twist $\sigma(G)$ introduced in \ref{my1}
is trivial (see \ref{sigmatr12}). This will imply that various maps used below are
$K$-orientable. 

We fix an orientation of the vector space $\Lie(G)$.

\subsubsection{}

We consider the correspondence of stacks 
\begin{equation}\label{cr12}[G\times G/G\times G]\stackrel{p}{\leftarrow}[G\times
G/G]\stackrel{q}{\rightarrow} [G/G]\ ,\end{equation}
where $p$ is induced by the identity on the level of spaces, and by the diagonal embedding of groups, and $q$ is given by the multiplication on the level of spaces, and by the identity on the level of groups.
We consider a twist  of $[G/G]$ of the form $\tau=\hol_*v^*\cT$.
By considering an equivalent correspondence of moduli spaces
associated to a pair of pants surface and \ref{ris4} we will obtain an isomorphism of
twists
$$r:q^*\tau\stackrel{\sim}{\rightarrow}p^*(\pr_1^*\tau+\pr_2^*\tau)\ .$$
We will furthermore construct a $K$-orientation of the proper and representable
map $q$ such that
\begin{equation}\label{prodde3}m:{}^{\tau} K([G/G])\otimes {}^{\tau} K([G/G])\rightarrow {}^{\tau}
K([G/G])\end{equation} defined by
$$m(x,y):=q_! r^*p^*(\pr_1^*x\cup \pr_2^*y)$$
is an associative unital product.
In fact, its complexification  will coincide with the product induced by the
identification
${}^{\tau} K([G/G])_\C\cong R(G)_\C/I$ given in \ref{rrg1}.

\subsubsection{}

Note that $[G/G]\rightarrow [*/G]$ is $-\sigma(G)$-$K$-orientable (see
\ref{my1}). By our assumption on $G$ we have
$\sigma(G)=0$.
Then $[G/G]\rightarrow [*/G]$ is $K$-orientable.
Therefore
$[G\times G/G\times G]\rightarrow [*/G\times G]$
 $K$-orientable.
By restriction to the diagonal subgroup we see that
$[G\times G/G]\rightarrow [*/G]$ is $K$-orientable.
It follows that $q$ is $K$-orientable.

\subsubsection{}\label{not51}

Let now $\Sigma$ be an oriented  pair of pants surface with ingoing boundary
components
$\partial_{i,\alpha} \Sigma$, $\alpha=1,2$, and outgoing boundary
component
$\partial_a \Sigma$. We assume that $\Sigma$ comes with a Riemannian
metric which has a product structure near the boundary such the
boundary circles are isometric to standard circles.
The correspondence
\begin{equation}\label{cr121}\cM(\partial_i\Sigma)\stackrel{p}{\leftarrow}\cM(\Sigma)\stackrel{q}{\rightarrow}
\cM(\partial_a\Sigma)\end{equation}
is equivalent to (\ref{cr12}).
Lemma \ref{ris4} together with \ref{sum41} now gives the desired isomorphism of twists
$$r:p^*\tau_i\stackrel{\sim}{\rightarrow} q^*\tau_a\ ,$$
where we use the notation of \ref{ia27}.
Note that this isomorphism may depend on the choice of the identification of
the correspondence (\ref{cr12}) with the correspondence (\ref{cr121}).

\subsubsection{}

We fix base points $b_\alpha\in \partial_{i,\alpha}\Sigma$,
$\alpha\in\{1,2\}$, and  $b\in \partial_a\Sigma$.
 Using the orientation of $\partial_i\Sigma$ opposite to
the induced one we define the holonomy map
$$\hol_i:\cM(\partial_i\Sigma)\rightarrow [G/G]\times [G/G]\ .$$
Let $\hol_a:\cM(\partial_a\Sigma)\rightarrow [G/G]$ be the holonomy
map associated to the outgoing boundary component.
The projection $[G/G]\rightarrow [*/G]$ induces a $R(G)$-module
structure on ${}^\tau K([G/G])$. Via the two projections $[G\times
G/G\times G]\rightarrow [G/G]\rightarrow [*/G]$ we have two
$R(G)$-module structures on ${}^{\pr_1^*\tau+\pr_2^*\tau}K([G\times
G/G\times G])$, which we write as left- and right actions.
Using the identifications via the holonomy maps we obtain
corresponding actions on ${}^{\tau_a}K(\cM(\partial_a\Sigma))$ and
${}^{\tau_i}K(\cM(\partial_i\Sigma))$.

\begin{lem}
The multiplication map $m$ is $R(G)$-bilinear.
\end{lem}
\proof
This is an immediate consequence of the commutativity of the diagram
$$\begin{array}{ccccc}
[G\times G/G\times G]&\stackrel{p}{\leftarrow}&[G\times
G/G]&\stackrel{q}{\rightarrow}& [G/G]\\
&\pr_\alpha\searrow&\downarrow&\swarrow&\\
&&{}[*/G]&&\end{array}
$$
and the projection formula (see \ref{asor5}).
\hB

\subsubsection{}

In the proof of Lemma \ref{rrg1} we have identified
${}^\tau K([G/G])_\C \cong R(G)_\C/I$ with the space of sections 
$\Gamma(S,V)$, where $V\rightarrow F^{reg}/W$ is a one-dimensional vector bundle associated to the character $\sign:W\rightarrow \{1,-1\}$, and $S\subset F^{reg}/W$ . Therefore
${}^{\pr_1^*\tau+\pr_2^*\tau} K([G\times G/G\times G])_\C\cong
\Gamma(S\times S,\pr_1^*V\otimes \pr_2^*V)$ with the left and right $\C[S]$-module structures
induced by the two projections $\pr_\alpha:F^{reg}\times F^{reg}\rightarrow F^{reg}$.
The multiplication thus induces a linear map
$$m_\C:\Gamma(S\times S,\pr_1^*V\otimes \pr_2^*V) \rightarrow \Gamma(S,V) \ .$$
Such a map is given by structure constants $C^{t}_{r,s}\in \Hom(V_r\otimes V_s,V_t)$, $r,s,t\in S$.
 
Since $m$ is $\C[S]$-bilinear we immediately conclude that
$C^t_{r,s}=0$ if not $t=s=r$.

We define the section  $c\in \Gamma(S,V^*)$ by
$c_s:=C^s_{s,s}$, $s\in S$, where $V^*$ is the dual bundle.

Note that all twisted $K$-groups are free $\Z$-modules and therefore
embed in their complexifications. In order to determine the product
(\ref{prodde3}) it
therefore suffices to calculate its complexification, i.e, the section $c$.
The following properties follow from the vanishing of the off-diagonal structure constants. 
\begin{kor}\label{first651}
The product (\ref{prodde3}) is associative and commutative.
\end{kor}
These properties in particular are independent on the choice of the $K$-orientation of $q$ and the choice of the isomorphism of twists $r$ (see \ref{not51} for the notation).

\subsubsection{}

In this subsection let us write $\Sigma_1$ for the pair of pants
surface considered above, $q_{1,a}$ for $q$, and $q_{1,i}$ for $p$. We furthermore consider an oriented surface
$\Sigma_0$ which is the union of a  disk and a cylinder such that it
has an ingoing boundary component (belonging to the cylinder) and two
outgoing boundary components.
We equip $\Sigma_0$ with Riemannian metric which has a product
structure such that the boundary components are isometric to the standard circle.
We fix an orientation reversing isometry
$\psi:\partial_{i}\Sigma_1\stackrel{\sim}{\rightarrow}
\partial_a\Sigma_0$ such that $\partial_{i,2}\Sigma_1$ is mapped to
the boundary of the disk.
This is exactly the situation considered in
\ref{nn11}. Let $\Sigma:=\Sigma_0\sharp_{\partial_a\Sigma_0\cong
  \partial_{i}\Sigma_1}\Sigma_1$ be the surface obtained by glueing.
Then $\Sigma$ is a cylinder.

The map
$q_{0,a}:\cM(\Sigma_0)\rightarrow \cM(\partial_a\Sigma_0)$ 
is equivalent to
$$[*/G]\times [G/G]\rightarrow [G/G]\times [G/G]\cong [G\times
G/G\times G]$$ and therefore proper, representable and $K$-orientable.
In fact, the choice of an orientation of $\Lie(G)$ induces a $K$-orientation of $q_{0,a}$.

Now observe that the  square in (\ref{se4318}) is cartesian.
Therefore the $K$-orientation of $q_{0,a}$ induces a $K$-orientation
of $j$, and together with the choice of a $K$-orientation of $q_{1,a}$
a $K$-orientation of $q_{a}:=q_{1,a}\circ j:\cM(\Sigma)\rightarrow
\cM(\partial_a\Sigma)$. 
Using Corollary \ref{rrc17} we obtain the  following identity of maps
${}^\tau K([G/G])\rightarrow {}^\tau K([G/G])$.
\begin{eqnarray}
m\circ (q_{0,a})_!\circ r_0^*\circ q_{0,i}^*&=&(q_{1,a})_!\circ
r_1^*\circ q_{1,i}^*\circ \psi^*\circ  (q_{0,a})_!\circ r_0^*\circ q_{0,i}^*\nonumber\\
&=&(q_{1,a})_!\circ
r_1^*\circ j_!\circ i^* \circ r_0^*\circ q_{0,i}^*\nonumber\\
&=&(q_{1,a})_!\circ j_! \circ j^*(r_1)^*\circ i^*(r_0)^* \circ i^* \circ q_{0,i}^*\nonumber\\
&=&(q_{a})_!\circ r^*\circ q_{i}^*\label{cl701}\ ,
\end{eqnarray}
where $r$ is associated to the cylinder.
Note that the correspondence (\ref{cne3}) associated to a cylinder
is equivalent to
$$[G/G]\stackrel{q_i}{\leftarrow}
  [G/G]\stackrel{q_a}{\rightarrow}[G/G]\ ,$$
where all maps are the identity.
In particular, there is a  distinguished $K$-orientation of
$q_{a}$, and a  distinguished  isomorphism of twists
$q_a^*\tau\stackrel{\sim}{\rightarrow} q_i^*\tau$.

 The set of $K$-orientations of $q_a$ (up to sign) is then identified with
$$H^2([G/G],\Z)\cong H^2(G\times_GEG,\Z)\subset  H^2(G,\Z)\cong
\Ext(\pi_1(G),\Z)$$
(Leray-Serre spectral sequence). 
Since $q_{1,a}\circ j$ is an isomorphism, 
 we see that
$j^*:H^2([G\times G/G],\Z)\rightarrow H^2([G/G],\Z)$ is surjective.
Therefore  we can choose the $K$-orientation of $q_{1,a}$ such that
the induced $K$-orientation
of $q_a=q_{1,a}\circ j$ is the distinguished one up to a sign.

Using again the surjectivity of $j^*$ we can adjust the isomorphism of twists
$r_1$ such that the induced isomorphism 
twists
$r:q_a^*\tau\stackrel{\sim}{\rightarrow} q_i^*\tau$ is  the
distinguished one.

We now fix the sign of the $K$-orientation of $q_{1,a}$ such that
the induced orientation of $q_a$ is the distinguished one.
In this case $(q_{a})_!\circ r^*\circ q_{i}^*=\id$. This fixes also
the class of $K$-orientations used to define the product $m$.

\subsubsection{}

Let $\Sigma^D\subset \Sigma_0$ be the component of the disk.
We indicate the related maps with the same superscript.
We have an element $E:=(q_a^D)_!\circ (r^D)^*(1)$.
The correspondence (\ref{cne3}) associated to the disk is 
equivalent to
$$*\leftarrow [*/G]\stackrel{q^D_a}{\rightarrow} [G/G]\ .$$
Therefore the element $E$ is the same as the one constructed
\ref{inr41}. The construction of $E$ in \ref{inr41} depends on  the choice of a section
$G\rightarrow \hat G$. Under the present assumptions on $G$ there is
only one such section, since $G$ has no non-trivial $U(1)$-valued characters.
In the present subsection we will see why we called $E$ the unit.

In fact the calculation \ref{cl701} and commutativity of the product \ref{first651} now gives
$$m(E,x)=m(x,E)=x\ , \quad \forall x\in {}^{\tau} K([G/G])\ .$$
We immediately conclude that the section $c$ is invertible and determined by
$$E\cong  c^{-1}\ ,$$ where we consider 
$E$ and $c^{-1}$ as sections of $V$.
The description of the product in terms of the basis
$(E_{[\chi]})_{[\chi]\in X_1(\hat T)/\hat W}$ is more complicated.

\subsubsection{}

We can state the final theorem about the product. We adopt the choices
of $K$-orientations fixed above. 
\begin{theorem}\label{pru78}
The product $m$ induces on ${}^\tau K([G/G])$ a commutative and
associative ring structure with identity $E$. 
Its complexification is isomorphic to the quotient $R(G)_\C/I$.
\end{theorem}

\end{document}